\documentclass[a4paper,11pt]{article}
\usepackage[latin1]{inputenc} 
\usepackage[T1]{fontenc}
\usepackage{lmodern}

\usepackage{amsthm,amsmath,amsfonts,amssymb,stmaryrd,bbm,mathrsfs}
\usepackage{mathtools}
\usepackage{enumitem}
\usepackage{url}
\usepackage[colorlinks=true, linkcolor=blue, urlcolor=black, citecolor=blue,pdfstartview=FitH]{hyperref}

\usepackage[french,english]{babel}
\usepackage{caption,tikz,color,subfigure}
\usetikzlibrary{shapes}
\usetikzlibrary{patterns}

\usepackage[top=2.7cm, bottom=2.7cm, left=2.5cm, right=2.5cm]{geometry}

\makeatletter
\renewcommand\theequation{\thesection.\arabic{equation}}
\@addtoreset{equation}{section}
\makeatother

\renewcommand{\captionfont}{\footnotesize}
\renewcommand{\captionlabelfont}{\footnotesize}

\setlist[enumerate,1]{label=(\roman*), font = \normalfont} 

\let\originalleft\left
\let\originalright\right
\renewcommand{\left}{\mathopen{}\mathclose\bgroup\originalleft}
\renewcommand{\right}{\aftergroup\egroup\originalright}

\newlength{\bibitemsep}
\newlength{\bibparskip}
\let\oldthebibliography\thebibliography
\renewcommand\thebibliography[1]{\oldthebibliography{#1}
  \setlength{\parskip}{\bibitemsep}
  \setlength{\itemsep}{\bibparskip}}

\newcommand{\N}{\mathbb{N}}
\newcommand{\Z}{\mathbb{Z}}
\newcommand{\Q}{\mathbb{Q}}
\newcommand{\R}{\mathbb{R}}
\newcommand{\C}{\mathbb{C}}
\newcommand{\T}{\mathbb{T}}
\renewcommand{\P}{\mathbb{P}}
\newcommand{\E}{\mathbb{E}}

\newcommand{\Ec}[1]{\mathbb{E} \left[#1\right]}
\newcommand{\Pp}[1]{\mathbb{P} \left(#1\right)}
\newcommand{\Ecsq}[2]{\mathbb{E} \left[#1\middle|#2\right]}
\newcommand{\Ppsq}[2]{\mathbb{P} \left(#1\middle|#2\right)}
\newcommand{\Varsq}[2]{\Var \left(#1\middle|#2\right)}

\newcommand{\1}{\mathbbm{1}}

\newcommand{\cN}{\mathcal{N}}
\newcommand{\cG}{\mathcal{G}}
\newcommand{\cB}{\mathcal{B}}
\newcommand{\cC}{\mathcal{C}}

\newcommand{\e}{\mathrm{e}}
\newcommand{\diff}{\mathop{}\mathopen{}\mathrm{d}}
\DeclareMathOperator{\Var}{Var}
\DeclareMathOperator{\Cov}{Cov}
\newcommand{\supp}{\mathrm{supp}}

\newcommand{\Exists}{\exists\,}
\newcommand{\Forall}{\forall\,}

\newcommand{\abs}[1]{\left\lvert#1\right\rvert}
\newcommand{\abso}[1]{\lvert#1\rvert}
\newcommand{\norme}[1]{\left\lVert#1\right\rVert}
\newcommand{\ps}[2]{\langle #1,#2 \rangle}

\newcommand{\petito}[1]{o\mathopen{}\left(#1\right)}
\newcommand{\grandO}[1]{O\mathopen{}\left(#1\right)}

\newcommand{\enstq}[2]{\left\{#1\mathrel{}\middle|\mathrel{}#2\right\}}
\newcommand{\restreinta}{\mathclose{}|\mathopen{}}

\newcommand{\relphantom}[1]{\mathrel{\phantom{#1}}}

\usepackage{anyfontsize}

\title{Two-temperatures overlap distribution for the 2D discrete Gaussian free field}

\author{
Michel \textsc{Pain}\thanks{DMA, \'Ecole Normale Sup\'erieure, PSL, CNRS, 75005 Paris, France 
\& Sorbonne Universit\'e, Sorbonne Paris Cit\'e, CNRS, Laboratoire de Probabilit\'es Statistique et Mod\'elisation, LPSM, F-75005 Paris, France. Email: \texttt{michel.pain@ens.fr}. Partially supported by ANR MALIN (ANR-16-CE93-0003).}
~and Olivier \textsc{Zindy}\thanks{Sorbonne Universit\'e, Sorbonne Paris Cit\'e, CNRS, Laboratoire de Probabilit\'es Statistique et Mod\'elisation, LPSM, F-75005 Paris, France. Email: \texttt{olivier.zindy@upmc.fr}. Partially supported by ANR MALIN (ANR-16-CE93-0003).}
}

\theoremstyle{plain}  
\newtheorem{thm}{Theorem}[section]
\newtheorem{prop}[thm]{Proposition}
\newtheorem{lem}[thm]{Lemma}

\theoremstyle{definition}

\theoremstyle{remark}
\newtheorem{rem}[thm]{Remark}

\begin{document}

\maketitle

\begin{abstract}
In this paper, we prove absence of temperature chaos for the two-dimensional discrete Gaussian free field using the convergence of the full extremal process, which has been obtained recently by Biskup and Louidor. 
This means that the overlap of two points chosen under Gibbs measures at different temperatures has a nontrivial distribution.
Whereas this distribution is the same as for the random energy model when the two points are sampled at the same temperature, we point out here that they are different when temperatures are distinct: more precisely, we prove that the mean overlap of two points chosen under Gibbs measures at different temperatures for the DGFF is strictly smaller than the REM's one.
Therefore, although neither of these models exhibits temperature chaos, one could say that the DGFF is more chaotic in temperature than the REM.
\end{abstract}

{
  \hypersetup{linkcolor=black}
  \tableofcontents
}

\section{Introduction}

\subsection{Spin glasses and chaos phenomenon}

The phenomenon of chaos in temperature or disorder is a classical problem in spin glasses, which was discovered by Fisher and Huse \cite{FisherHuse1986} for the Edwards-Anderson model, and Bray and Moore~\cite{BrayMoore1987} for the Sherrington-Kirkpatrick model. It arose from the discovery that, for some models, a small change in the external parameters (such as temperature or disorder) can induce a dramatic change of the overall energy landscape or may modify the location of the ground state and the organization of the pure states for the Gibbs measure.
In the last decades, this phenomenon has been studied extensively by physicists for various models, see Rizzo \cite{Rizzo2009} for a recent survey. 
Several mathematical results have also been obtained by Chen and Panchenko~\cite{ChenPanchenko2013}, Chen \cite{Chen2014}, Auffinger and Chen \cite{AuffingerChen2016}, Panchenko \cite{Panchenko2016}, Subag \cite{subag2017}, Chen and Panchenko \cite{ChenPanchenko2017} and Ben Arous, Subag and Zeitouni \cite{BenArousSubagZeitouni2018arxiv} in recent years.

This paper concerns the problem of chaos in temperature, which can be described more precisely as follows. 
The \textit{overlap} between two configurations is defined as the normalized covariance between energies of these configurations.
\textit{Temperature chaos} means that, if one samples independently two spin configurations from Gibbs measures at different temperatures but with fixed disorder, then their overlap will be almost deterministic. 
This phenomenon happens when the set of most likely configurations under the Gibbs measure changes radically when the temperature varies only slightly: in spin glass models, the Gibbs measure at a given temperature is concentrated near some fixed level of energy and on some pure states 
and both of them can change with temperature. 

Some spin glass models display temperature chaos and some others do not.
For example, Subag \cite{subag2017} proved recently the absence of chaos in temperature for spherical {\it pure} $p$-spin models with $p \geq 3$, while Ben Arous, Subag and Zeitouni \cite{BenArousSubagZeitouni2018arxiv} proved chaos in temperature for some spherical {\it mixed} $p$-spin models. 
Both results rely on a precise geometrical description of the Gibbs measure and hold at low temperature ($\beta$ large enough), in part of the one-step replica symmetry breaking (1-RSB) regime. 
For the spherical pure $p$-spin models, the supports of the Gibbs measures are close to each other for different (low enough) temperatures, while, for spherical mixed $p$-spin models, the Gibbs measure concentrates on thin spherical bands which depend on the temperature. 
This difference explains partially why chaos in temperature occurs for the second class of models but not for the first one.

In order to shed some light on the mysteries of the Parisi theory for mean field spin glasses, Derrida introduced in the 80's the random energy model (REM) \cite{Derrida1981}, where the Gaussian energy levels are assumed to be independent, and its generalization, the generalized random energy model (GREM) \cite{Derrida1985}, whose correlations are given by a tree structure of finite depth. These tractable models have been extensively studied and allowed, in particular, to investigate the phenomenon of replica symmetry breaking. We refer to Bolthausen \cite{Bolthausen2007} and Kistler \cite{Kistler2015} for connection to spin glass theory. Let us mention here that Kurkova \cite{kurkova2003} proved the absence of chaos in temperature for the REM; this will be made more precise later.

Natural hierarchical models with an infinite number of levels are the branching Brownian motion (BBM) and the branching random walk (BRW), see e.g. the seminal paper by Derrida and Spohn \cite{derridaspohn88}, who introduced directed polymers on trees (a BRW with i.i.d.\@ displacements) as an infinite hierarchical extension of the GREM for spin glasses. Recently this field has experienced a revival with many remarkable contributions and repercussions in other areas (cover times \cite{BeliusKistler2017,BeliusRosenZeitouni2020,schmidt2020}, characteristic polynomials of random unitary matrices \cite{ArguinBeliusBourgade2017,ChhaibiMadauleNajnudel2018,LambertPaquette2016arxiv,PaquetteZeitouni2018}, the Riemann zeta function on the critical line \cite{Najnudel2018,abbrs2018} and a random model of the Riemann zeta function \cite{ArguinBeliusHarper2017,ArguinTai2019,Ouimet2018}). We refer to Shi \cite{shi2015} for a survey on BRW and to Bovier \cite{Bovier2015} for a source motivated by connections between BBM and spin glasses.

Physicists suggested that Gaussian BRW and BBM should belong to a universality class called log-correlated Gaussian fields. We refer to the works by Carpentier et Le Doussal \cite{carpentier-ledoussal2001}, Fyodorov and Bouchaud \cite{FyodorovBouchaud2008a,FyodorovBouchaud2008b} and Fyodorov, Le Doussal and Rosso \cite{FyodorovLeDoussalRosso2009} for connection with spin glass theory.
Such a model is the two-dimensional discrete Gaussian free field that we study in this paper and will be described precisely in the next subsection. It appears that this model has an implicit hierachical structure similar to BRW. The lecture notes of Biskup \cite{biskup2020} give a general and excellent account of recent results about two-dimensional discrete Gaussian free field. See also \cite{zeitouni2017notes} for connections with BRW.

In this paper, we prove absence of temperature chaos for the two-dimensional discrete Gaussian free field using the convergence of the full extremal process recently proved by Biskup and Louidor \cite{biskuplouidor2018}. We also show that the mean overlap between two points sampled independently according to Gibbs measures at different temperatures is strictly smaller than the REM's one, which might be surprising since the overlap distribution is the same for both models if the two points are sampled at the same temperature, see \eqref{eq:arguin_zindy} and \eqref{eq:same_distribution_overlap}.

\subsection{The two-dimensional discrete Gaussian free field}

We consider in this paper the two-dimensional discrete Gaussian free field (DGFF) on an \textit{admissible} lattice approximation of a bounded open set $D \subset \C$.
More precisely, let $D$ be a bounded open set of $\C$ such that its boundary $\partial D$ has only a finite number of connected components, each of which has a positive diameter and is a $\cC^1$ path (i.e.\@ the range of a $\cC^1$ function from $[0,1]$ to  $\C$).
We denote by $\mathrm{dist}_\infty$ the $\ell^\infty$-distance on $\Z^2$.
Let $(D_N)_{N \geq 1}$ be the sequence of subsets of $\Z^2$ defined by
\[
D_N \coloneqq \left\{ x \in \Z^2 : \mathrm{dist}_\infty \left( \frac{x}{N},D^c \right) > \frac{1}{N} \right\}.
\]
These assumptions are slightly more restrictive than those in the recent papers by Biskup and Louidor \cite{biskuplouidor2020, biskuplouidor2019, biskuplouidor2018}, but only to avoid some technical issues due to boundary effects. 
The \textit{discrete Gaussian free field} on $D_N$ is the Gaussian process $(h^N_x)_{x\in D_N}$ with covariance given by the Green function $G_N$ of the simple random walk on $\Z^2$ killed upon exiting $D_N$. See \eqref{eq:Greenfunction} for a detailed definition of $G_N$.
In the sequel, we will skip the dependence in $N$ for the field $h^N$, denoting simply $(h_x)_{x\in D_N}$.
In comparison with spin glass theory, $D_N$ plays the role of the configurations set and $-h_x$ is the energy of configuration $x$.
Moreover, the \textit{overlap} between $x,y \in D_N$ can be defined by
\[
q_N(x,y) \coloneqq \frac{\E[h_x h_y]}{\max_{z \in D_N} \E[h_z^2]}
= \frac{G_N(x,y)}{\max_{z \in D_N} G_N(z,z)}.
\]
Note that $q_N(x,y) \in [0,1]$ and behaves roughly like $1-\frac{\log \Vert x-y \Vert}{\log N}$ at least for points $x$ and $y$ far enough from the boundary $\partial D_N$, using the asymptotic behavior of $G_N$ (see Lemma \ref{lemma:overlap-distance}).
It was shown by Bolthausen, Deuschel and Giacomin \cite{bdg2001} that the maximum of the DGFF, or {\it ground state} in spin glass terminology, satisfies
\[
\frac{\max_{x\in V_N} h_x}{\log N^2} 
\xrightarrow[N\to\infty]{} \sqrt{g}, \quad \text{in probability,}
\]
with $g \coloneqq 2/\pi$ and in the special case of the square lattice $V_N := (0,N)^2 \cap \Z^2$.
Their technique was later refined by Daviaud \cite{Daviaud2006} who computed the {\it log-number of high points} in $V_N$: for any $0<\lambda<1$,
\begin{align} \label{eq:number_of_high_points}
\lim_{N\to\infty} \frac{1}{\log N^2}  \log \#\{x \in V_N: h_x \geq  \lambda \sqrt{g} \log N^2\}= 1-\lambda^2, \quad \text{in probability.}
\end{align}
This result holds also with a general admissible lattice approximation $(D_N)_{N\geq 1}$, as a consequence of \cite[Theorem 2.1]{biskuplouidor2019}, which sharpens considerably the results of~\cite{Daviaud2006}.

\begin{figure}[t]
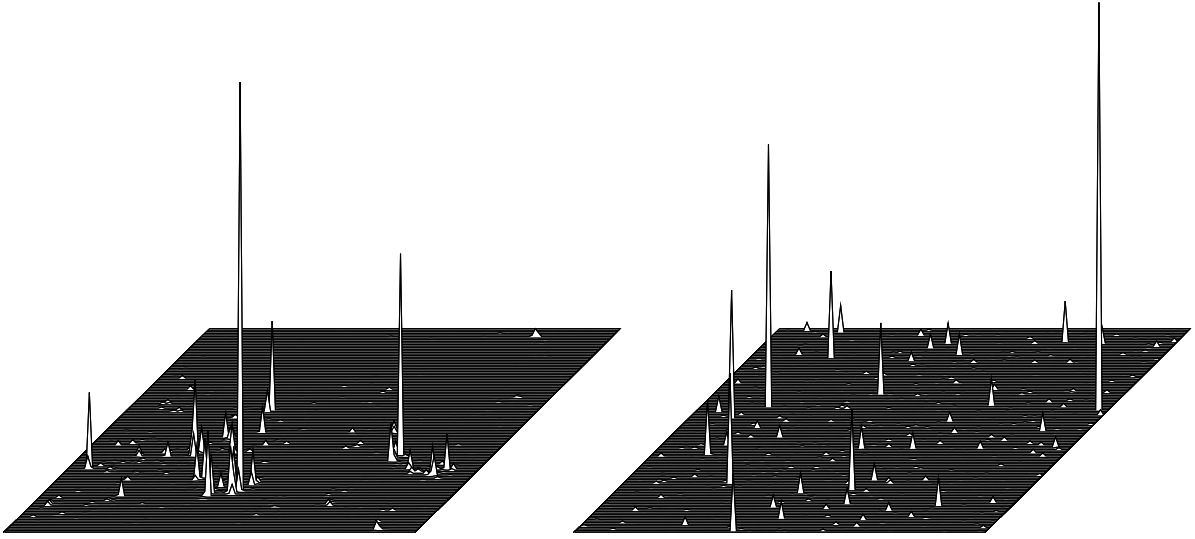

\vspace{-0.5cm}
\centering
\subfigure{
 
}
\caption{Realization of $\e^{\beta h}$, where $h$ is a DGFF on the left and a REM on the right. The domain is the square $V_N = (0,N) \cap \Z^2$ with $N = 127$ and $\beta$ is slightly supercritical ($\beta = 1.1 \cdot \beta_c$). Note that points carrying the mass of the Gibbs measure at inverse temperature $\beta$ are gathered in clusters for the DGFF and are uniformly distributed on the square for the REM.}
\label{fig:exp_GFF_REM}
\end{figure}

The \textit{Gibbs measure} at inverse temperature $\beta > 0$ associated with the DGFF is defined by
\[
\cG_{\beta,N} 
\coloneqq \frac{1}{Z_{\beta,N}} \sum_{x \in D_N} \e^{\beta h_x} \delta_x ,
\]
where $Z_{\beta,N} \coloneqq \sum_{x \in D_N} \e^{\beta h_x}$ is the partition function.
This model exhibits a phase transition in the asymptotic behavior of the \textit{free energy} defined by 
\[
f_N(\beta) \coloneqq  \frac{\log Z_{\beta,N}}{\log N^2}.
\]
Indeed, it follows easily from \eqref{eq:number_of_high_points} (see \cite{arguinzindy2015}) that the free energy converges as $N \to \infty$ and
\begin{align} \label{eq:limit_free_energy}
\lim_{N\to\infty} f_N(\beta) 
= f(\beta) \coloneqq  \left\{
\begin{array}{ll}
1 + (\beta/\beta_c)^2, & \text{if } \beta \leq \beta_c, \\
2 \beta/\beta_c, & \text{if } \beta \geq \beta_c,
\end{array}
\right.
\quad \text{in probability and in } L^1,
\end{align}
where the critical inverse temperature is given by 
\[
\beta_c \coloneqq \sqrt{2\pi} = \frac{2}{\sqrt{g}}.
\]
Arguin and Zindy \cite{arguinzindy2015} used the previous results to prove that the model displays a one-step replica symmetry breaking in spin glass terminology. More precisely, they showed that, on the square lattice $V_N$,
\begin{align} \label{eq:arguin_zindy}
\E \left[ \cG_{\beta,N}^{\otimes 2} (q_N(u,v) \in \cdot) \right] 
\xrightarrow[N\to\infty]{} 
\frac{\beta_c}{\beta}\delta_0 + \left(1-\frac{\beta_c}{\beta} \right)\delta_1,
\quad \Forall \, \beta > \beta_c,
\end{align}
where, on the left-hand side, $u$ and $v$ are chosen independently according to $\cG_{\beta,N}$ given the Gaussian free field $h$, before taking the expectation w.r.t.\@ $h$.
In other words, the left-hand side is the pushforward of the measure $\E[ \cG_{\beta,N}^{\otimes 2}]$ on $D_N^2$ by the function $q_N$.
To prove this, they adopted a spin glass approach, introducing a perturbed DGFF and adapting the so-called Bovier--Kurkova technique \cite{bovierkurkova2004-1,bovierkurkova2004-2}, which relates the free energy of the perturbed DGFF with the overlap distribution of the initial DGFF. 
Let us notice that all the aforementioned results hold also for the REM, which we define, for each $N \geq 1$, as a family $(h^\mathrm{REM}_x)_{x \in D_N}$ of i.i.d.\@ centered Gaussian random variables with variance $\max_{x\in D_N} G_N(x,x)$, in order to be comparable to the DGFF. 
A different behavior between the REM and the DGFF appears in the study of the second order of the maximum: 
it has been proved by Bolthausen, Deuschel and Zeitouni \cite{bdz2011}, Bramson and Zeitouni \cite{bramsonzeitouni2012} and Bramson, Ding and Zeitouni \cite{bdz2016-2}, also in the special case of the square lattice $V_N$, that the centered maximum of the DGFF $\max_{x \in V_N} h_x - m_N$, where
\begin{align} \label{eq:cv_max_dgff}
m_N \coloneqq 2 \sqrt{g} \log N - \frac{3}{4} \sqrt{g} \log \log N,
\end{align}
converges in distribution towards a nontrivial limit, whereas the right centering term for the REM is $m_N^\mathrm{REM} \coloneqq 2 \sqrt{g} \log N - \frac{1}{4} \sqrt{g} \log \log N$.
Biskup and Louidor described the limit in \cite{biskuplouidor2016} and extended the result to a class of {\it admissible} lattice approximations for $D$ in \cite{biskuplouidor2020}.
The next step in the description of extremes of the DGFF has been the study of the extremal process, which describes the field seen from position $m_N$.
More precisely, the \textit{full extremal process} is defined as the following random measure on $D \times \R \times \R^{\Z^2}$:
\begin{align} \label{eq:def_extremal_point_process}
\eta_{N,r} 
\coloneqq \sum_{x \in D_N}
\1_{\{h_x = \max_{y\in\Lambda_r(x)} h_y\}}
\delta_{x/N} \otimes \delta_{h_x-m_N} 
\otimes \delta_{(h_x - h_{x+z})_{z \in \Z^2}}, 
\end{align}
where $\Lambda_r(x) \coloneqq \{ y \in D_N : \norme{x-y} \leq r \}$.
It encodes the rescaled position of the local maxima, their centered value and the field seen from this local maximum.
The convergence of the full extremal process has been proved by Biskup and Louidor \cite{biskuplouidor2018} and can be stated as follows.
There exist a random finite Borel measure $Z^D$ on $\overline{D}$ and a probability measure $\nu$ on $(\R_+)^{\Z^2}$ such that, for any sequence $(r_N)_{N\geq1}$ of positive real numbers with $r_N \to \infty$ and $N/r_N \to \infty$,
\begin{align} \label{eq:biskup_louidor_full_extremal_process_convergence}
\eta_{N,r_N} 
\xrightarrow[N\to\infty]{\text{law}} 
\eta \coloneqq \mathrm{PPP} \left( Z^D(\diff z) \otimes \e^{-\beta_c h} \diff h
\otimes \nu(\diff \phi) \right),
\end{align}
in the sense that, for any continuous function $f \colon \overline{D} \times (\R \cup \{\infty\}) \times \overline{\R}^{\Z^2} \to \R$ with compact support, $\eta_{N,r_N}(f) = \int f \diff \eta_{N,r_N}$ converges in distribution towards $\eta(f)$, where we set $\overline{\R} \coloneqq \R \cup \{-\infty,\infty\}$.
It means that, given $Z^D$, rescaled positions of local maxima are asymptotically i.i.d.\@ variables with law $Z^D/Z^D(\overline{D})$, their centered values are given by an independent Poisson point process with intensity $Z^D(\overline{D}) \e^{-\beta_c h} \diff h$ and to each of these local maxima is attached an independent cluster with law $\nu$.
The random measure $Z^D$ has been studied in \cite{biskuplouidor2020} and identified as the critical Liouville quantum gravity obtained by exponentiating the continuous GFF on $D$: in particular, $Z^D$ has full support in $\overline{D}$ and $Z^D(\partial D) = 0$ almost surely.
Moreover, Biskup and Louidor \cite{biskuplouidor2018} give an explicit description of the cluster law~$\nu$, see Subsection \ref{subsection:decorations_are_nontrivial} for more details.
The convergence of the full extremal process leads to a precise description of the supercritical Gibbs measure $\cG_{\beta,N}$ for $\beta > \beta_c$.
For large $N$, this measure is concentrated on points $x \in D_N$ such that $h_x = m_N + O(1)$.
Therefore, in spin glass terminology, the pure states are balls of diameter $O(1)$ centered at each local maximum of height $m_N+O(1)$: if two points are in the same pure state, then they have an overlap close to 1 and, if they are in two different pure states, their overlap is close to 0.
This explains the one-step replica symmetry breaking of the model.
Moreover, the ordered weights of the pure states under $\cG_{\beta,N}$ follows asymptotically a \textit{Poisson--Dirichlet distribution} of parameter $(\beta_c/\beta)$ denoted $\mathrm{PD}(\beta_c/\beta)$, which is the law of $(a_i/\sum_j a_j, i\geq 1)_\downarrow$ where $\downarrow$ stands for the decreasing rearrangement and $(a_i)_{i\geq 1}$ are the atoms of a Poisson point process on $(0,\infty)$ with intensity  $s^{-(\beta_c/\beta)-1} \diff s$.
More precisely, Biskup and Louidor \cite[Corollary 2.7]{biskuplouidor2018} proved that, for $\beta>\beta_c$, on the space of Radon measures on $\overline{D}$ endowed with the vague convergence, we have
\begin{align} \label{eq:Gibbs-PD}
\sum_{x \in D_N} \cG_{\beta,N}(\{x\}) \, \delta_{x/N}
\xrightarrow[N\to\infty]{} \sum_{i\geq 1} p_i \, \delta_{\chi_i},
\quad \text{in distribution},
\end{align}
where, conditionally on $Z^D$, $(\chi_i)_{i\geq 1}$ is an i.i.d.\@ sequence of random variables with distribution $Z^D/Z^D(D)$ and $(p_i)_{i\geq 1}$ is independent of $(\chi_i)_{i\geq 1}$ with law $\mathrm{PD}(\beta_c/\beta)$.
We also mention that a weaker version of this Poisson--Dirichlet convergence (for the overlap distribution) has been obtained by Arguin and Zindy \cite{arguinzindy2014,arguinzindy2015}.
For the REM, the supercritical Gibbs measure $\cG_{\beta,N}^\mathrm{REM}$ for $\beta > \beta_c$ is also carried by extremal points, which are the $x \in D_N$ such that $h_x^\mathrm{REM} = m_N^\mathrm{REM} + O(1)$. 
Moreover, these points are uniformly chosen in $D_N$ and, in the limit, their heights are given by an independent Poisson point process with intensity $c \, \e^{-\beta_c h} \diff h$, with $c > 0$.
Hence, the pure states are here the singletons formed by extremal points and are at distance of order $N$ of each other.
The pure states have also Gibbs weights following asymptotically the Poisson--Dirichlet distribution of parameter $\beta_c/\beta$ and therefore the overlap under $(\cG_{\beta,N}^\mathrm{REM})^{\otimes 2}$ has asymptotically the same law as for the DGFF. 
See Figure \ref{fig:exp_GFF_REM}.
Note that results in that direction have also been obtained by Bovier and Kurkova \cite{bovierkurkova2004-1,bovierkurkova2004-2} for more elaborate variants of the REM: the GREM and the CREM. 

\subsection{Overlap distribution at two different temperatures}

In the previous subsection, we saw that the pure states of the supercritical Gibbs measure of the DGFF are points with height $m_N + O(1)$. In particular, they do not depend on the temperature and, for $\beta,\beta' > \beta_c$, two points chosen independently according to $\cG_{\beta,N}$ and $\cG_{\beta',N}$ have a positive probability to be in the same pure state, but also to be in different pure states: therefore, their overlap can be either $0$ or $1$ and there is clearly no temperature chaos for the DGFF.
This is made more precise in our first result, stating the convergence of the overlap under $\cG_{\beta,N} \otimes \cG_{\beta',N}$.
\begin{thm} \label{theorem:convergence_overlap_GFF}
Let $\beta, \beta' > 0$.
\begin{enumerate}
\item If $\beta \leq \beta_c$ or $\beta' \leq \beta_c$, then, for all $a \in (0,1)$,
\[
\cG_{\beta,N} \otimes \cG_{\beta',N} (q_N(u,v) \geq a)
\xrightarrow[N\to\infty]{} 0, \quad  \text{in $L^1$}.
\]
\item  If $\beta > \beta_c$ and $\beta'  >\beta_c$, then, for all $a \in (0,1)$,
\[
\cG_{\beta,N} \otimes \cG_{\beta',N} (q_N(u,v) \geq a)
\xrightarrow[N\to\infty]{} 
Q(\beta,\beta'), \quad \text{in distribution},
\]
where
\[
Q(\beta,\beta') \coloneqq
\frac{\sum_{k\geq1} 
\left( \sum_{x \in \Z^2} \e^{\beta (\xi_k-\phi^k_x)} \right)
\left( \sum_{x \in \Z^2} \e^{\beta' (\xi_k-\phi^k_x)} \right)}
{\left( \sum_{k\geq1} \sum_{x \in \Z^2} \e^{\beta (\xi_k-\phi^k_x)} \right)
\left( \sum_{k\geq1} \sum_{x \in \Z^2} \e^{\beta' (\xi_k-\phi^k_x)} \right)},
\]
with $(\xi_k)_{k\geq1}$ the atoms of a Poisson point process with intensity $\e^{-\beta_c h} \diff h$ and $(\phi^k)_{k\geq1}$ independent samples from the measure $\nu$, independent of $(\xi_k)_{k\geq1}$.
\end{enumerate}
\end{thm}
In other words, this result proves the convergence of the pushforward of the measure $\cG_{\beta,N} \otimes \cG_{\beta',N}$ on $D_N^2$ by the function $q_N$, which is a random measure on $[0,1]$. The limit is either $\delta_0$ if $\beta \wedge \beta' \leq \beta_c$, or $(1-Q(\beta,\beta')) \delta_0 + Q(\beta,\beta') \delta_1$ otherwise.
In case (ii), given the limit of the full extremal process, $Q(\beta,\beta')$ is simply the probability of choosing two points in the same cluster, when they are chosen proportionally to their Gibbs weights with inverse temperature $\beta$ and $\beta'$ respectively.
For the REM, the overlap between $x,y \in D_N$ is defined by
\[
q_N^\mathrm{REM}(x,y) \coloneqq \frac{\E \left[h_x^\mathrm{REM} h_y^\mathrm{REM}\right]}{\max_{z \in D_N} G_N(z,z)} 
= \1_{\{x=y \}}.
\]
Kurkova \cite{kurkova2003} proved the following result: if $\beta,\beta' > \beta_c$, then, for any $a \in (0,1)$, we have
\begin{equation} \label{eq:REM}
\cG_{\beta,N}^\mathrm{REM} \otimes \cG_{\beta',N}^\mathrm{REM} (q_N^\mathrm{REM}(u,v) \geq a)
\xrightarrow[N\to\infty]{\mathrm{(d)}} 
Q^\mathrm{REM}(\beta,\beta')
\coloneqq
\frac{\sum_{k\geq 1} \e^{\beta \xi_k} \e^{\beta' \xi_k}
}{\left( \sum_{k\geq 1} \e^{\beta \xi_k} \right)
	\left( \sum_{k\geq 1} \e^{\beta' \xi_k} \right)},
\end{equation}
where the $(\xi_k)_{k\geq 1}$ are also the atoms of a Poisson point process with intensity $\e^{-\beta_c x} \diff x$.
Our aim is now to compare $Q(\beta,\beta')$ and $Q^\mathrm{REM}(\beta,\beta')$.
In the case $\beta = \beta'$, it is known that $Q(\beta,\beta)$ and $Q^\mathrm{REM}(\beta,\beta)$ have the same distribution.
Indeed, for $\beta > \beta_c$ and $k\geq1$, we introduce the random variable (well-defined by \cite[Theorem 2.6]{biskuplouidor2018})
\begin{align*}
X_{\beta,k} \coloneqq \frac{1}{\beta} \log \sum_{x \in \Z^2} \e^{-\beta \phi^k_x},
\end{align*}
so that the limit of the overlap can be rewritten as
\begin{align} \label{eq:reformulation_Q}
Q(\beta,\beta') 
= \frac{\sum_{k\geq1} \e^{\beta (\xi_k+X_{\beta,k})} 
\e^{\beta' (\xi_k+X_{\beta',k})}}
{\left( \sum_{k\geq1} \e^{\beta (\xi_k+X_{\beta,k})} \right)
\left( \sum_{k\geq1} \e^{\beta' (\xi_k+X_{\beta',k})} \right)}.
\end{align}
Then, noting that $(X_{\beta,k})_{k\geq1}$ is a sequence of i.i.d.\@ random variables, which is independent of $(\xi_k)_{k\geq1}$, and that $\E[\e^{\beta_c X_{\beta,k}}] < \infty$ by \cite[Theorem 2.6]{biskuplouidor2018}, we have 
\begin{align} \label{eq:chgt_PPP}
\left(\xi_k + X_{\beta,k} \right)_{k\geq 1}
\overset{\text{(d)}}{=}
\left( \xi_k + \beta_c^{-1} \log \E \left[ \e^{\beta_c X_{\beta,1}} \right] \right)_{k\geq 1},
\end{align}
where the equality in distribution holds with both sides seen as point processes (see for example \cite[Proposition 8.7]{bolthausensznitman2002}).
It shows that
\begin{align} \label{eq:same_distribution_overlap}
Q(\beta,\beta) 
\overset{\text{(d)}}{=}
\frac{\sum_{k\geq1} \e^{2\beta \xi_k} }
{\left( \sum_{k\geq1} \e^{\beta \xi_k} \right)^2
}
= Q^\mathrm{REM}(\beta,\beta).
\end{align}
Moreover, it follows from a simple change of variable that $(\e^{\beta \xi_i}/ \sum_{k\geq1} \e^{\beta \xi_k},i\geq 1)_\downarrow$ has the same law as $(p_i)_{i\geq 1}$ a Poisson--Dirichlet random variable of parameter $(\beta_c/\beta)$, 
hence $Q(\beta,\beta)$ and $Q^\mathrm{REM}(\beta,\beta)$ have the same law as $\sum_{i\geq 1} p_i^2$, which corresponds to the fact that the pure states have Poisson--Dirichlet weights.

One may ask whether $Q(\beta,\beta')$ and $Q^\mathrm{REM}(\beta,\beta')$ have the same distribution when $\beta \neq \beta'$.
The answer is negative and our second result shows that, in mean, the overlap is less likely to be close to 1 under $\cG_{\beta,N} \otimes \cG_{\beta',N}$ that under $\cG_{\beta,N}^\mathrm{REM} \otimes \cG_{\beta',N}^\mathrm{REM}$.
\begin{thm} \label{theorem:comparaison_overlap_GFF}
For any $\beta, \beta' > \beta_c$ such that $\beta \neq \beta'$, we have
\[
\Ec{Q(\beta,\beta')}
< \Ec{Q^\mathrm{REM}(\beta,\beta')}.
\]
\end{thm}
The reason behind this inequality is the following. 
For the REM, the weight of a pure state depends only on its height: a likely pure state under $\cG_{\beta,N}^\mathrm{REM}$ will also be likely under $\cG_{\beta',N}^\mathrm{REM}$.
Rather, the weight of a pure state for the DGFF depends both on the height of the local maximum and on the geometry of the cluster around it, hence some pure states are likely for $\beta$ close to $\beta_c$ and unlikely for large $\beta$, or vice-versa. 
Therefore, it is more difficult to choose twice the same pure state under the Gibbs measures at two different temperatures for the DGFF than for the REM and one could say that the DGFF is more chaotic in temperature than the REM. 
This result leads to the picture displayed in Figure \ref{fig:curves} for the mean overlap in both models.
\begin{figure}[ht]
\centering
\begin{tikzpicture}[scale = 2] 
\draw[->,>=latex] (0,0) -- (5.1,0) node[right]{$\beta$};
\draw[->,>=latex] (0,0) -- (0,1.5*1.1);
\draw (-0.025,1.5) -- (0.025,1.5);
\draw (0,1.5) node[left]{$1$};
\draw (0,0) node[below left]{$0$};


\draw[thick,color=blue!70] (0,0) -- (1,0);	
\draw[thick,color=red!90,dashed] (0,0) -- (1,0);	

\draw[thick,color=blue!70] plot[domain=1:2,samples=100] 
	(\x,{1.5*sqrt(\x-1)/2});
\draw[thick,color=red!90] plot[domain=1.00001:2,samples=100] 
	(\x,{1.5*9/20*exp(1/3*ln(\x - 1)) + 1.5*1/20*exp(2*ln(\x - 1))});
	
\draw[thick,color=blue!70] plot[domain=2:5,samples=100] 
	(\x,{1.5*7/10-1.5*1/5*exp(15/4*ln(3/(\x+1)))});	
\draw[thick,color=red!90] plot[domain=2:5,samples=100] 
	(\x,{1.5*3/4-1.5*1/4*exp(6*ln(6/(\x+4)))});

\draw[color=blue!70,->,>=latex] (3.3,0.93) -- (3.55,0.75);
\draw[color=blue!70] (3.5,0.8) node[below right]{$\E[Q(\beta,\beta')]$};
\draw[color=red!90,->,>=latex] (3.15,1.05) -- (2.9,1.22);	
\draw[color=red!90] (2.95,1.17) node[above left]{$\E[Q^\text{REM}(\beta,\beta')]$};	


\draw[fill] (1,0) circle (0.02);
\draw (1,0) node[below]{$\beta_c$};

\draw[dashed] (2,-0.025) -- (2,1.5*0.5);
\draw (2,0) node[below]{$\beta'$};
\draw[dashed] (-0.025,1.5*0.5) -- (2,1.5*0.5);
\draw (0,1.5*0.5) node[left]{$\beta_c/\beta'$};
\draw[fill] (2,1.5*0.5) circle (0.02);

\end{tikzpicture}
\caption{Schematic representation of the mean overlap at two different inverse temperatures $\beta$ and $\beta'$, as $\beta'$ is fixed and $\beta$ varies, for the DGFF (in blue) and the REM (in red). The difference between the curves is actually exaggerated on the figure: computer calculations suggest that they should be indistinguishable at this scale.}
\label{fig:curves}
\end{figure}
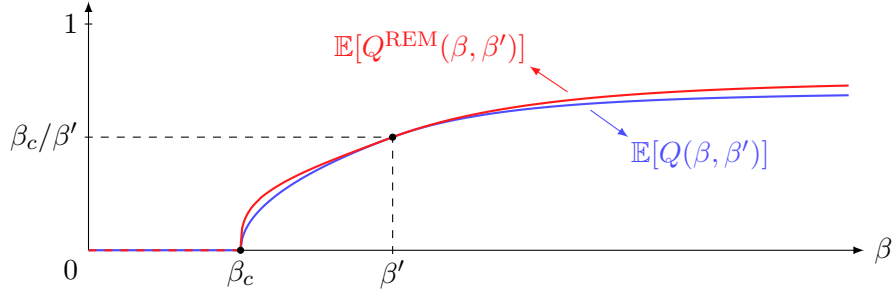
\begin{rem} 
Theorem \ref{theorem:comparaison_overlap_GFF} can be seen as the comparison of $\E[\ps{p}{q}] = \E[\sum_{i\geq 1} p_i q_i]$, for two different couples $(p,q)$ with marginal laws $\mathrm{PD}(\beta_c/\beta)$ and $\mathrm{PD}(\beta_c/\beta')$.
In the case of the REM, the couple is characterized by the relation 
$q_i = p_i^{\beta'/\beta}/ (\sum_{k\geq 1} p_k^{\beta'/\beta})$, and one could ask if this particular choice reaches the minimum in the Wasserstein distance between the measures $\mathrm{PD}(\beta_c/\beta)$ and $\mathrm{PD}(\beta_c/\beta')$ on the Hilbert space $\ell^2$ of square-summable sequences.
The answer is negative by \cite[Theorem 2.3]{CuestaMatran1989} and the fact that 
$\P \otimes \P \left( \{ (\omega,\omega') : 
\ps{p(\omega)-p(\omega')}{q(\omega)-q(\omega')} < 0 
\} \right) > 0$ in that case.
\end{rem}
\begin{rem} \label{rem:beta_infinite}
In the case where $\beta < \infty$ and $\beta' = \infty$, it is easier to see that the overlap distribution for the DGFF should be different from the REM's one.
We can define $\cG_{\infty,N}$ as the measure which gives mass $1$ to the point where $\max_{x\in D_N} h_x$ is reached. Then, assuming that $(\xi_k)_{k\geq 1}$ is ranked in the decreasing order, one can prove in this case that the limiting overlap distribution is given by 
\[
Q(\beta,\infty) = \frac{\sum_{x \in \Z^2} \e^{\beta (\xi_1-\phi^1_x)}}{\sum_{k\geq1} \sum_{x \in \Z^2} \e^{\beta (\xi_k-\phi^k_x)} }
= \frac{\e^{\beta (\xi_1+X_{\beta,1})} }
{\sum_{k\geq1} \e^{\beta (\xi_k+X_{\beta,k})}},
\]
which is asymptotically the probability that the point chosen according to $\cG_{\beta,N}$ is in the same cluster as the highest particle.
On the other hand, we have 
$Q^\mathrm{REM} (\beta,\infty) = \e^{\beta \xi_1} / (\sum_{k\geq1} \e^{\beta \xi_k})$.
Then, it follows from \eqref{eq:chgt_PPP} that $Q(\beta,\infty)$ is stochastically strictly dominated by $Q^\mathrm{REM} (\beta,\infty)$, because with positive probability $\max_{k\geq 1} (\xi_k+X_{\beta,k})$ is not reached at $k = 1$ (this is true as soon as $X_{\beta,k}$ is not a.s.\@ constant, which holds by the proof of Lemma \ref{lem:decoration_non_triviales}).
\end{rem}
\begin{rem}
Similar results could be proved for the branching Brownian motion or the branching random walk. 
Indeed, using the convergence of the extremal process of the BBM \cite{abbs2013,abk2013,bovierhartung2017} or of the BRW \cite{madaule2017,Mallein2018}, one can deduce Theorem \ref{theorem:convergence_overlap_GFF}.(ii) (see in particular the proof of Theorem 4.3 in \cite{Mallein2018}).
For Theorem \ref{theorem:comparaison_overlap_GFF}, the method used here is quite general and works also for BBM and BRW, aside from Lemma \ref{lem:decoration_non_triviales} that should be adapted.
\end{rem}

\subsection{Organization of the paper}

The paper is organized as follows. Theorem \ref{theorem:convergence_overlap_GFF} is proved in Section \ref{section:convergence_overlap} and Theorem \ref{theorem:comparaison_overlap_GFF} in Section \ref{section:comparaison}. Appendix \ref{appendix} contains well-known results concerning the two-dimensional Gaussian free field.

\section{Convergence of the overlap distribution}
\label{section:convergence_overlap}

\subsection{Proof of Part (i) of Theorem \ref{theorem:convergence_overlap_GFF}}

The proof of Part (i) of Theorem \ref{theorem:convergence_overlap_GFF} will be a direct consequence of the following result, which relies on the fact that the free energy contains all information about the mean overlap under the measure $\E[\cG_{\beta,N}^{\otimes 2}]$.
\begin{prop} \label{prop:overlap-hautetemperature}
If $\beta \leq \beta_c$, then for any $a \in (0,1]$,
 \[
\cG_{\beta,N}^{\otimes 2} (q_N(u,v) \geq a)
\xrightarrow[N\to\infty]{} 0, \quad \text{in } L^1.
\]
\end{prop}
\begin{proof}
Recall that $f_N(\beta)$ denotes the free energy of the DGFF on $D_N$ at inverse temperature~$\beta$.
Observe that $\E[f_N(\beta)]$ is a convex function of the parameter $\beta.$ 
Moreover, it was already mentionned that $\E [ f_N(\beta)]$ converges towards $f(\beta)$, which is also convex in $\beta$. 
In particular, by a standard result of convexity (see e.g. Proposition I.3.2 in \cite{simon1993}), at each point of differentiability of $f$, the limit 
of the derivatives equals the derivative of the limit. 
It follows from \eqref{eq:limit_free_energy} that $f$ is differentiable at any $\beta>0$, hence we get
\begin{align*}
f'(\beta)
= \lim_{N\to\infty} \frac{\diff}{\diff \beta} \E \left[ f_N(\beta) \right]
= \lim_{N\to\infty} \frac{1}{\log N^2} \sum_{x \in D_N} 
	\Ec{\frac{h_x \e^{\beta h_x}}{\sum_{z \in D_N} \e^{\beta h_z}}}.
\end{align*}
Applying Gaussian integration by part (see Lemma \ref{lem:Gaussian-IbP}) with respect to the factor $h_x$, it follows that
\begin{align*} 
f'(\beta)
= \lim_{N\to\infty} \frac{\beta}{\log N^2} \left(
\sum_{x \in D_N} \Ec{h_x^2} 
	\Ec{\frac{\e^{\beta h_x}}{\sum_{z \in D_N} \e^{\beta h_z}}}
- \sum_{x,y \in D_N} \Ec{h_x h_y} 
	\Ec{\frac{\e^{\beta h_x}\e^{\beta h_y}}{(\sum_{z \in D_N} \e^{\beta h_z})^2}}
\right).
\end{align*}
In order to deal with the first sum, we introduce $D_{N,\delta} \coloneqq \{ x \in D_N : d(x,D_N^c) > N^{1-\delta} \}$ for some small $\delta \in (0,1)$.
It follows from \eqref{eq:relation_G_N_and_a} and \eqref{eq:asymptotics_a}, that $\max_{x \in D_N} G_N(x,x) = \frac{2}{\pi} \log N + O_N(1)$ and $G_N(x,x) \geq \frac{2 (1- \delta)}{\pi} \log N + O_N(1)$ uniformly in $x \in D_{N,\delta}$ as $N \to \infty$.
On the other hand, we have $\E[\cG_{\beta,N}(D_{N,\delta}^c)] \to 0$ as $N\to\infty$ by Lemma \ref{lemma:complement-negligible}.
Combining this and letting $\delta \to 0$, we get
\[
f'(\beta) = \lim_{N\to\infty} \frac{\beta}{\pi} \left(1-\E \left[ \cG_{\beta,N}^{\otimes 2}[q_N(u,v)]\right]\right).
\]
Moreover, it follows from \eqref{eq:limit_free_energy} that $f'(\beta)=\beta/\pi$ for any $\beta \leq \beta_c,$ and therefore
\[
\lim_{N\to\infty} \E \left[ \cG_{\beta,N}^{\otimes 2}[q_N(u,v)]\right]=0,
\]
which concludes the proof.
\end{proof}
\begin{proof}[Proof of Part (i) of Theorem \ref{theorem:convergence_overlap_GFF}]
By symmetry, we can assume that $\beta \le \beta_c$. 
Moreover, by Lemma \ref{lemma:overlap-distance}.(i), it is enough to prove that $\cG_{\beta,N} \otimes \cG_{\beta',N} (\lVert u-v \rVert \leq N^{1-a}) 
\to 0$ as $N \to \infty$ in $L^1(\P)$, for any $a \in (0,1]$.
For this, we introduce a partition of $D_N$ by considering the intersection of $D_N$ with a partition of $\R^2$ with boxes of side-length $N^{1-a}$: let $\cB$ denote this partition of $D_N$.
Then, we have the following upper bound
\begin{align*}
& \cG_{\beta,N} \otimes \cG_{\beta',N} (\lVert u-v \rVert \leq N^{1-a}) \\
& \leq \sum_{B \in \cB} \cG_{\beta,N}(B) 
\cG_{\beta',N}(\{ y \in D_N : \Exists x \in B : \Vert x-y \Vert \leq N^{1-a} \})
  \\
& \leq \left( \max_{B \in \cB} \cG_{\beta,N}(B) \right)
\sum_{B \in \cB} 
\cG_{\beta',N}(\{ y \in D_N : \Exists x \in B : \Vert x-y \Vert \leq N^{1-a} \})
  \\
& \leq 9 \, \left( \max_{B \in \cB} \cG_{\beta,N}(B) \right) \sum_{B \in \cB} \cG_{\beta',N}(B) 
= 9 \max_{B \in \cB} \cG_{\beta,N}(B),
\end{align*}
noting that the factor $9$, which appears in the last inequality, comes from the fact that if $x \in B$ and $\Vert x-y \Vert \leq N^{1-a}$, then $y$ belongs either to $B$ or to one of the $8$ boxes neighbor to $B$,
see Figure \ref{figure:boxes}.
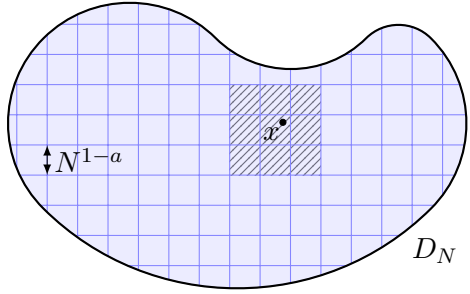
\begin{figure}[t]
\centering
\begin{tikzpicture}[scale = 1]
\fill [pattern=north east lines] (2.8,0.8) -- (4,0.8) -- (4,2) -- (2.8,2) -- cycle;
\begin{scope}
\clip[rotate=-45,shift={(0,0.1)}] (1.4,5) arc (360:270:1.4) -- (0,3.6) arc (90:270:1.6) -- (0,0.4) arc (270:360:3.6) -- (3.6,4) arc (0:90:1.6) -- (2,5.6) arc (90:180:0.6);
\draw[blue!80,step=0.4] (-2,-2) grid (6,4);
\end{scope}
\draw (5.5,-0.2) node{$D_N$};
\draw[rotate=-45,shift={(0,0.1)},fill=blue!15,opacity=0.5] (1.4,5) arc (360:270:1.4) -- (0,3.6) arc (90:270:1.6) -- (0,0.4) arc (270:360:3.6) -- (3.6,4) arc (0:90:1.6) -- (2,5.6) arc (90:180:0.6);
\draw[rotate=-45,shift={(0,0.1)},thick] (1.4,5) arc (360:270:1.4) -- (0,3.6) arc (90:270:1.6) -- (0,0.4) arc (270:360:3.6) -- (3.6,4) arc (0:90:1.6) -- (2,5.6) arc (90:180:0.6);
\draw[fill] (3.5,1.5) circle (0.04);
\draw (3.35,1.35) node{$x$}; 
\draw[<->,>=latex] (0.4,0.8) -- (0.4,1.2);
\draw (0.95,1) node{$N^{1-a}$};
\end{tikzpicture}
\caption{The set $D_N$ is partitionned using a square grid of span $N^{1-a}$. 
If $x$ is in a given box and $\Vert x-y \Vert \leq N^{1-a}$, then $y$ belongs to one of the hatched boxes.}
\label{figure:boxes}
\end{figure}
On the other hand, for any $B \in \cB$, we have
\begin{align*}
\left( \max_{B \in \cB} \cG_{\beta,N}(B) \right)^2
= \max_{B \in \cB} \cG_{\beta,N}^{\otimes 2} (B^2)
\leq \cG_{\beta,N}^{\otimes 2} ( \lVert u-v \rVert \leq \sqrt{2} N^{1-a})
\xrightarrow[N\to\infty]{L^1(\P)} 0,
\end{align*}
where the convergence follows from Proposition \ref{prop:overlap-hautetemperature} combined with Lemmas \ref{lemma:overlap-distance}.(ii) and \ref{lemma:complement-negligible}.
This concludes the proof.
\end{proof}

\subsection{Proof of Part (ii) of Theorem \ref{theorem:convergence_overlap_GFF}}

For proving Part (ii) of Theorem \ref{theorem:convergence_overlap_GFF}, we use the convergence of the full extremal process \eqref{eq:biskup_louidor_full_extremal_process_convergence}.
We first prove finite-dimensional convergence of the following family of random measures on~$D$:
\begin{align*}
\rho_{\beta,N} 
\coloneqq \sum_{x \in D_N} \e^{\beta(h_x-m_N)} \delta_{x/N},
\quad \Forall \beta > \beta_c.
\end{align*}
The one-dimensional convergence has been proven by Biskup and Louidor \cite[Theorem 2.6]{biskuplouidor2018}: on the space of Radon measures on $\overline{D}$ endowed with the vague convergence, we have
\begin{align} \label{eq:theorem2.6_biskup_louidor}
\rho_{\beta,N} 
\xrightarrow[N\to\infty]{} 
c(\beta) Z^D(D)^{\beta/\beta_c} 
\Sigma_\beta,
\quad \text{in distribution},
\end{align}
where $\Sigma_\beta$ is a random measure on $\overline{D}$ such that, conditionally on $Z^D$, $\Sigma_\beta = \sum_{k\geq1} \e^{\beta \xi_k} \delta_{\chi_k}$, with $(\chi_k)_{k\geq1}$ an i.i.d.\@ sequence of random variables with distribution $Z^D/Z^D(D)$, $(\xi_k)_{k\geq1}$ a Poisson point process with intensity $\e^{-\beta_c h} \diff h$ independent of $(\chi_k)_{k\geq1}$ and $c(\beta)$ a positive constant defined by
$c(\beta) \coloneqq \beta^{-\beta/\beta_c} \
(\E[( \sum_{x\in \Z^2} \e^{-\beta \phi_x})^{\beta_c/\beta}])^{\beta/\beta_c}$ where $\phi$ has distribution $\nu$.
We introduce another family of random measures $(\rho_\beta)_{\beta > \beta_c}$ defined as follows.
Conditionally on $Z^D$, let $(\chi_k)_{k\geq1}$ be an i.i.d.\@ sequence of random variables with distribution $Z^D/Z^D(D)$, $(\xi_k)_{k\geq1}$ a Poisson point process with intensity $\e^{-\beta_c h} \diff h$ and $(\phi^k)_{k\geq1}$ be an i.i.d.\@ sequence of random variables with distribution $\nu$ such that $(\chi_k)_{k\geq1}$, $(\xi_k)_{k\geq1}$ and $(\phi^k)_{k\geq1}$ are independent.
Then, we set 
\begin{align*}
\rho_\beta 
\coloneqq Z^D(D)^{\beta/\beta_c} \sum_{k\geq1} \e^{\beta \xi_k} 
\left( \sum_{x \in \Z^2} \e^{-\beta \phi^k_x} \right) 
\delta_{\chi_k}, 
\quad  \Forall \beta > \beta_c.
\end{align*}
By \eqref{eq:chgt_PPP}, $(c(\beta) Z^D(D)^{\beta/\beta_c} \Sigma_\beta)_{\beta>\beta_c}$ and $(\rho_\beta)_{\beta > \beta_c}$ have the same one-dimensional marginal distributions, but not the same multi-dimensional ones.
Our first aim is to prove the following extension of convergence \eqref{eq:theorem2.6_biskup_louidor} to the case of several temperatures.
\begin{prop} \label{prop:cv_gff_mesure_de_gibbs_spatiale}
For any $p \geq 1$ and $\beta_1, \dots, \beta_p > \beta_c$, we have
\begin{align*}
\rho_{\beta_1,N} \otimes \dots \otimes \rho_{\beta_p,N} 
\xrightarrow[N\to\infty]{} 
\rho_{\beta_1} \otimes \dots \otimes \rho_{\beta_p} ,
\quad \text{in distribution},
\end{align*}
on the space of Radon measures on $\overline{D}^p$ endowed with the vague convergence.
\end{prop}
Before proving this proposition, we first recall two results from \cite{biskuplouidor2018}:
by \cite[Equation (6.58)]{biskuplouidor2018}, for any $\beta > \beta_c$ and $\eta > 0$,
\begin{align} \label{eq:height_in_partition_function}
\lim_{\ell \to\infty} \limsup_{N\to\infty} 
\P \left( \sum_{x\in D_N} \e^{\beta (h_x-m_N)} \1_{\{h_x - m_N \notin [-\ell,\ell]\}} > \eta \right) = 0
\end{align}
and by \cite[Lemma B.11]{biskuplouidor2018}, for any $K>0$,
\begin{align} \label{eq:distance_extremal_points}
\lim_{r\to\infty} \limsup_{N\to\infty} \P \left( \Exists x,y \in D_N : h_x,h_y > m_N -K, r < \norme{x-y} < N/r \right) = 0.
\end{align}
We can now proceed to the proof of Proposition \ref{prop:cv_gff_mesure_de_gibbs_spatiale}.
\begin{proof}[Proof of Proposition \ref{prop:cv_gff_mesure_de_gibbs_spatiale}]
The goal is to apply properly the convergence of the full extremal process \eqref{eq:biskup_louidor_full_extremal_process_convergence} and, for this, we use arguments very close to those in the proof of Theorem~2.6 of~\cite{biskuplouidor2018}.
In order to prove the result, it is sufficient to show that, for any $f_1,\dots, f_p \colon \overline{D} \to \R_+$ continuous functions with compact support (which is always the case since $\overline{D}$ is compact),
\begin{align*}
\left( \rho_{\beta_1,N}(f_1), \dots, \rho_{\beta_p,N}(f_p) \right)
\xrightarrow[N\to\infty]{} 
\left( \rho_{\beta_1}(f_1), \dots, \rho_{\beta_p}(f_p) \right),
\quad \text{in distribution},
\end{align*}
see Kallenberg \cite[Lemma 4.1]{kallenberg2017}. 
For this, it is sufficient to prove that
\begin{align} \label{ba}
\Ec{\exp \left( - \sum_{i=1}^p \rho_{\beta_i,N}(f_i) \right)}
\xrightarrow[N\to\infty]{} 
\Ec{\exp \left( - \sum_{i=1}^p \rho_{\beta_i}(f_i) \right)}.
\end{align}
Proceeding as in \cite{biskuplouidor2018}, in order to apply \eqref{eq:biskup_louidor_full_extremal_process_convergence}, we need to keep in $\rho_{\beta,N}$ only the $x\in D_N$ such that $\abs{h_x - m_N}$ is bounded.
Therefore, we introduce an increasing sequence of continuous functions $(g_\ell)_{\ell \geq 1}$ from $\R$ to $\R$ such that, for each $\ell \geq 1$, $\1_{[-\ell,\ell]} \leq g_\ell \leq \1_{[-\ell-1,\ell+1]}$.
It follows from \eqref{eq:height_in_partition_function} that, for any $\varepsilon >0$, $\beta > \beta_c$ and $f \colon \overline{D} \to \R_+$ continuous function, we have
\begin{align} \label{bh}
\lim_{\ell\to\infty} \limsup_{N\to\infty} 
\Pp{ \abs{\rho_{\beta,N}(f) 
	- \sum_{x \in D_N} \e^{\beta(h_x-m_N)} f(x/N) g_{\ell-1} (h_x-m_N) }
> \varepsilon }
= 0.
\end{align}
Then, we introduce the event 
\[
A_{N,\ell,r} \coloneqq \{ \Forall x,y \in D_N, (h_x,h_y > m_N -\ell) \Rightarrow \Vert x-y \Vert \notin (r, N/r) \}.
\]
We are going to transform repeatedly the sum in \eqref{bh}, with changes that are well-controlled on the event $A_{N,\ell,r}$. By \eqref{eq:distance_extremal_points}, $A_{N,\ell,r}$ is very likely as $N \to \infty$ and then $r \to \infty$, so this will be sufficient.

Let $(r_N)_{N\geq1}$ be a sequence of positive real numbers with $r_N \to \infty$ and $N/r_N \to \infty$, that will be used to apply \eqref{eq:biskup_louidor_full_extremal_process_convergence}.
On the event $A_{N,\ell,r}$, for $N$ large enough such that $r \leq r_N < \frac{N}{2r}$, we have
\begin{align}
& \sum_{x \in D_N} \e^{\beta(h_x-m_N)} f(x/N) g_{\ell-1} (h_x-m_N) \nonumber \\
& = \sum_{x \in D_N}
\1_{\{h_x = \max_{y\in\Lambda_{r_N}(x)} h_y\}} \e^{\beta(h_x-m_N)}
\sum_{y\in\Lambda_r(x)} \e^{\beta(h_y-h_x)} f(y/N) g_{\ell-1}(h_y-m_N), \label{bg}
\end{align}
where we emphasize that $\Lambda_{r_N}(x)$ and $\Lambda_r(x)$ do not have the same index.
Note that we know $\sum_{x \in D_N} \e^{\beta(h_x-m_N)}$ is tight by \cite[Equation (6.57)]{biskuplouidor2018} combined with the tightness of $\max_{x \in V_N} h_x - m_N$.
Hence, we can replace $f(y/N)$ with $f(x/N)$ in the right-hand side of \eqref{bg}: the error committed on $A_{N,\ell,r}$ is smaller than $\omega_f(r/N) \sum_{x \in D_N} \e^{\beta(h_x-m_N)}$, where $\omega_f$ is a modulus of continuity for $f$ (which is uniformly continuous), and hence tends to $0$ in probability as $N \to \infty$.
Besides this change, note that we can substitute $g_{\ell-1}(h_y-m_N)$ for $g_{\ell}(h_x-m_N) g_{2\ell}(h_y-h_x)$ so that \eqref{bg} is replaced with
\begin{align*}
\sum_{x \in D_N}
\1_{\{h_x = \max_{y\in\Lambda_{r_N}(x)} h_y\}} \e^{\beta(h_x-m_N)}
\sum_{y\in\Lambda_r(x)} \e^{\beta(h_y-h_x)} f(x/N) g_{\ell}(h_x-m_N) g_{2\ell}(h_y-h_x).
\end{align*}
Indeed, on the event $A_{N,\ell,r} \cap \{ \max_{x \in D_N} h_x - m_N \leq \ell \}$, we have, using that $h_y \leq h_x$, 
\begin{align*}
g_{\ell-1}(h_y-m_N)
& \leq \1_{\{h_y-m_N \in [-\ell,\ell] \}}
\leq \1_{\{h_x-m_N \in [-\ell,\ell] \}} 
	\1_{\{h_x-h_y \in [-2\ell,2\ell] \}} \\
& \leq g_{\ell}(h_x-m_N) g_{2\ell}(h_y-h_x),
\end{align*}
and, therefore, on this likely event (by \eqref{eq:cv_max_dgff} and \eqref{eq:distance_extremal_points}), the error term can be bounded by
$\norme{f}_\infty \sum_{x \in D_N} \e^{\beta(h_x-m_N)} \1_{\{h_y-m_N \notin [-\ell+1,\ell-1] \}}$ which tends to 0 in probability as $N \to \infty$ and $\ell \to \infty$ by \eqref{eq:height_in_partition_function}.
Finally, we proved the following analogue of \cite[Equation (6.59)]{biskuplouidor2018}:
\begin{align} \label{bb}
\lim_{\ell\to\infty} \limsup_{r\to\infty} \limsup_{N\to\infty} 
\Pp{ \abs{\rho_{\beta,N}(f) - \ps{\eta_{N,r_N}}{F_{\beta,r,\ell}} }
> \varepsilon }
= 0,
\end{align}
where we set $F_{\beta,r,\ell} (x,h,\phi) \coloneqq \e^{\beta h} f(x) g_{\ell}(h) Y_{r,\ell}^\beta(\phi)$
and $Y_{r,\ell}^\beta (\phi) \coloneqq \sum_{y \in \Lambda_r(0)} \e^{-\beta \phi_y} g_{2 \ell}(\phi_y)$ in order to get
\begin{align*}
\ps{\eta_{N,r_N}}{F_{\beta,r,\ell}}
& = \sum_{x \in D_N}
\1_{\{h_x = \max_{y\in\Lambda_{r_N}(x)} h_y\}} \e^{\beta(h_x-m_N)}
f(x/N) g_{\ell}(h_x-m_N)
Y_{r,\ell}^\beta \left( h_x - h_{x + \cdot} \right),
\end{align*}
recalling that $\eta_{N,r_N}$ is the extremal point process defined in \eqref{eq:def_extremal_point_process}.
Now, we apply what precedes and the full extremal process convergence in order to prove \eqref{ba}.
Let $F^i_{r,\ell} (x,h,\phi) \coloneqq \e^{\beta_i h} f_i(x) g_{\ell}(h) Y_{r,\ell}^{\beta_i}(\phi)$ for each $i \in \{ 1,\dots,p \}$.
It follows from \eqref{bb} that
\begin{align*}
\lim_{\ell\to\infty} \limsup_{r\to\infty} \limsup_{N\to\infty} 
\abs{ \Ec{\exp \left( - \sum_{i=1}^p \rho_{\beta_i,N}(f_i) \right)} 
- \Ec{\exp \left( - \sum_{i=1}^p \ps{\eta_{N,r_N}}{F^i_{r,\ell}} \right)} }
= 0.
\end{align*}
But, applying the full extremal process convergence \eqref{eq:biskup_louidor_full_extremal_process_convergence} to the function $\sum_{i=1}^p F^i_{r,\ell}$ which is continuous with compact support from $\overline{D} \times \R \times \R^{\Lambda_r(0)} \to \R_+$, we have
\begin{align}
& \Ec{\exp \left( - \sum_{i=1}^p \ps{\eta_{N,r_N}}{F^i_{r,\ell}} \right)} 
\nonumber \\
& \xrightarrow[N\to\infty]{}
\Ec{\exp \left( - \int_{\overline{D} \times \R \times \R^{\Z^2}}
Z^D(\diff x) \otimes \e^{-\beta_c h} \diff h \otimes \nu(\diff \phi)
\left( 1- \e^{-\sum_{i=1}^p F^i_{r,\ell}(x,h,\phi)} \right)
\right)}.\label{bc}
\end{align}
Note that $F^i_{r,\ell}(x,h,\phi) \uparrow F^i(x,h,\phi) \coloneqq \e^{\beta_i h} f_i(x) Y^{\beta_i}(\phi)$ as $r \to \infty$ and then $\ell \to \infty$,
where $Y^{\beta}(\phi) \coloneqq \sum_{y \in \Z^2} \e^{-\beta \phi_y}$.
Therefore, taking the limit inside the expectation (by dominated convergence theorem) and then inside the integral (by monotone convergence theorem), yields that the right-hand side of \eqref{bc} tends to 
\begin{align}
\Ec{\exp \left( - \int_{\overline{D} \times \R \times \R^{\Z^2}}
Z^D(\diff x) \otimes \e^{-\beta_c h} \diff h \otimes \nu(\diff \phi)
\left( 1-\e^{-\sum_{i=1}^p F^i(x,h,\phi)} \right)
\right)}, \label{bd}
\end{align}
as $r \to \infty$ and then $\ell \to \infty$.
Finally, note that \eqref{bd} is equal to the right-hand side of \eqref{ba} and, therefore, the result is proved.
\end{proof}
Before proceeding to the proof of Part (ii) of Theorem \ref{theorem:convergence_overlap_GFF}, note that it follows from \eqref{eq:theorem2.6_biskup_louidor} that the renormalized partition function, for $\beta > \beta_c$, converges in distribution towards a positive limit:
\begin{align} \label{eq:convergence_partition_function}
\e^{-\beta m_N} Z_{\beta,N} = \sum_{x \in D_N} \e^{\beta (h_x-m_N)}
\xrightarrow[N \to \infty]{} 
c(\beta) Z^D(D)^{\beta/\beta_c} \Sigma_\beta(\overline{D}),
\quad \text{in distribution}.
\end{align}
Combining this with \eqref{eq:height_in_partition_function}, it shows that the supercritical Gibbs measure is mainly supported by extremal points: more precisely, for any $\beta > \beta_c$ and $\eta > 0$,
\begin{align} \label{eq:height_in_Gibbs}
\lim_{\ell\to\infty} \limsup_{N\to\infty} 
\P( \cG_{\beta,N}( h_u - m_N \notin [-\ell,\ell]) > \eta)
= 0.
\end{align}
Finally, combining this with \eqref{eq:distance_extremal_points}, it follows that, for any $\beta,\beta' > \beta_c$, 
\begin{align} \label{eq:distance_under_supercritical_Gibbs}
\lim_{r\to\infty} \limsup_{N\to\infty} 
\P( \cG_{\beta,N} \otimes \cG_{\beta',N} (r < \norme{u-v} < N/r) > \eta)
= 0,
\end{align}
which means that two points sampled accordingly to two supercritical Gibbs measures are typically either very close to each other or very far.
\begin{proof}[Proof of Part (ii) of Theorem \ref{theorem:convergence_overlap_GFF}]
Our first aim is to prove that, for any $a \in (0,1)$ and $\eta > 0$,
\begin{align}
\limsup_{\varepsilon \to 0}
\limsup_{N\to\infty}
\Pp{
\abs{ \cG_{\beta,N} \otimes \cG_{\beta',N} (q_N(u,v) \geq a)
- \cG_{\beta,N} \otimes \cG_{\beta',N} (\Vert u-v \Vert \leq \varepsilon N) }
> \eta }
= 0. \label{be}
\end{align}
Fix some $\delta \in (0,1-a)$ and recall that $D_{N,\delta} = \{ x \in D_N : d(x,D_N^c) > N^{1-\delta} \}$.
It follows from Lemma \ref{lemma:overlap-distance} that, for any $\varepsilon > 0$, we have, for $N$ large enough,
\begin{align*}
\{ \Vert x-y \Vert \leq \varepsilon^{-1} \} \cap D_{N,\delta}^2
\subset
\{ q_N(x,y) \geq a \} \cap D_{N,\delta}^2
\subset
\{ \Vert x-y \Vert \leq \varepsilon N \} \cap D_{N,\delta}^2
\end{align*}
and it follows that
\begin{align}
& \abs{ \cG_{\beta,N} \otimes \cG_{\beta',N} (q_N(u,v) \geq a)
- \cG_{\beta,N} \otimes \cG_{\beta',N} (\Vert u-v \Vert \leq \varepsilon N) } \nonumber \\
& \leq \cG_{\beta,N} \otimes \cG_{\beta',N} (D_N^2 \setminus D_{N,\delta}^2)
+ \cG_{\beta,N} \otimes \cG_{\beta',N} (\varepsilon^{-1} <\Vert u-v \Vert  \leq \varepsilon N). \label{bf}
\end{align}
The first term in the right-hand side of \eqref{bf} is smaller than $\cG_{\beta,N} (D_{N,\delta}^c) + \cG_{\beta',N} (D_{N,\delta}^c)$, which tends to 0 in $\P$-probability as $N \to \infty$ by Lemma \ref{lemma:complement-negligible}.
On the other hand, the second term in the right-hand side of \eqref{bf} tends to zero in probability as $N \to \infty$ and then $\varepsilon \to 0$, by \eqref{eq:distance_under_supercritical_Gibbs}.
Therefore, \eqref{be} is proved.
For any  $\varepsilon > 0$, let $f_\varepsilon \colon D^2 \to \R_+$ be a continuous function such that $\1_{\{\abs{z-z'} \leq \varepsilon\}} \leq f_\varepsilon(z,z') \leq \1_{\{\abs{z-z'} \leq 2\varepsilon\}}$.
Then, we have 
\begin{align*}
\frac{\rho_{\beta,N} \otimes \rho_{\beta',N}(f_{\varepsilon/2})}{\rho_{\beta,N} \otimes \rho_{\beta',N}(D^2)}
\leq \cG_{\beta,N} \otimes \cG_{\beta',N} (\Vert u-v \Vert  \leq \varepsilon N) 
\leq \frac{\rho_{\beta,N} \otimes \rho_{\beta',N}(f_\varepsilon)}{\rho_{\beta,N} \otimes \rho_{\beta',N}(D^2)}.
\end{align*}
Moreover, by Proposition \ref{prop:cv_gff_mesure_de_gibbs_spatiale}, we get
\begin{align*}
\frac{\rho_{\beta,N} \otimes \rho_{\beta',N}(f_\varepsilon)}{\rho_{\beta,N} \otimes \rho_{\beta',N}(D^2)}
\xrightarrow[N\to\infty]{\text{(d)}} 
\frac{\rho_{\beta} \otimes \rho_{\beta'}(f_\varepsilon)}{\rho_{\beta} \otimes \rho_{\beta'}(D^2)}
\xrightarrow[\varepsilon\to0]{\text{a.s.}} 
\frac{\rho_{\beta} \otimes \rho_{\beta'} (\{ (z,z) : z \in D \})}
{\rho_{\beta} \otimes \rho_{\beta'} (D^2)}.
\end{align*}
Coming back to \eqref{be}, this implies that 
\begin{align*}
\cG_{\beta,N} \otimes \cG_{\beta',N} (q_N(u,v) \geq a)
\xrightarrow[N\to\infty]{\text{(d)}} 
\frac{\rho_{\beta} \otimes \rho_{\beta'} (\{ (z,z) : z \in D \})}
{\rho_{\beta} \otimes \rho_{\beta'} (D^2)},
\end{align*}
which is equal to $Q(\beta,\beta')$, concluding the proof of Part (ii) of Theorem \ref{theorem:convergence_overlap_GFF}.
\end{proof}

\section{Comparison with the overlap for the REM}
\label{section:comparaison}

\subsection{Structure of the proof of Theorem \ref{theorem:comparaison_overlap_GFF}}

In this subsection, we split the proof of Theorem \ref{theorem:comparaison_overlap_GFF} into three lemmas that will be proved in the next subsections.
Therefore, we consider some $\beta, \beta' > \beta_c$ such that $\beta \neq \beta'$.
Let $\phi$ be a random field with distribution $\nu$ and, for any $\beta > \beta_c$, $X_{\beta} \coloneqq \frac{1}{\beta} \log \sum_{x \in \Z^2} \e^{-\beta \phi_x}$.
Recall from \eqref{eq:reformulation_Q} that
\[
Q(\beta,\beta') 
= \frac{\sum_{k\geq1} \e^{\beta (\xi_k+X_{\beta,k})} 
\e^{\beta' (\xi_k+X_{\beta',k})}}
{\left( \sum_{k\geq1} \e^{\beta (\xi_k+X_{\beta,k})} \right)
\left( \sum_{k\geq1} \e^{\beta' (\xi_k+X_{\beta',k})} \right)},
\]
with $(\xi_k)_{k\geq1}$ the atoms of a Poisson point process with intensity $\e^{-\beta_c h} \diff h$ and $(X_{\beta,k},X_{\beta',k})_{k\geq1}$ i.i.d.\@ copies of $(X_{\beta},X_{\beta'})$, independent of $(\xi_k)_{k\geq1}$. 
In the case $\beta = \beta'$, we saw in the introduction that the classical argument to get rid of the $X_{\beta,k}$'s in the description of the overlap distribution is the fact that $(\xi_k + X_{\beta,k})_{k\geq 1}$ has the same distribution as 
$(\xi_k + c_\beta)_{k\geq 1}$ where $c_\beta$ is a constant (see \eqref{eq:chgt_PPP}).
Here, we have instead to work with the joint distribution $(\xi_k + X_{\beta,k},\xi_k + X_{\beta',k})_{k\geq 1}$ for $\beta \neq \beta'$. 
The following lemma shows that, if we apply the previous change of point process to the first coordinate, the random shift in the second coordinate is still independent of $(\xi_k)_{k\geq 1}$, although this can seem counter-intuitive at first sight.
This result appeared previously in a paper by Panchenko and Talagrand \cite[Lemma 2.1]{panchenkotalagrand2007-1}.
It will be proved for the sake of completeness in Subsection \ref{subsection:change_of_point_process}. 

\begin{lem} \label{lem:changement_de_processus_ponctuel}
Let $c_\beta \coloneqq \beta_c^{-1} \log \E[\e^{\beta_c X_{\beta}}]$ and
$Y$ be a random variable such that, for any measurable function $f \colon \R \to \R_+$,
\begin{align} \label{eq:def_Y}
\Ec{f(Y)}  = \frac{\Ec{\e^{\beta_c X_{\beta}} f(X_{\beta'} - X_{\beta})}}{\Ec{\e^{\beta_c X_{\beta}}}}.
\end{align}
Then, the point process $(\xi_k + X_{\beta,k},\xi_k + X_{\beta',k})_{k\geq 1}$ with values in $\R^2$ has the same distribution as $(\xi_k + c_\beta,\xi_k + c_\beta + Y_k)_{k\geq 1}$, where $(Y_k)_{k\geq 1}$ are i.i.d.\@ copies of $Y$ independent of $(\xi_k)_{k\geq 1}$.
\end{lem}

Now, our aim is to prove that the random perturbations $Y_k$ on the second coordinate of $(\xi_k + c_\beta,\xi_k + c_\beta + Y_k)_{k\geq 1}$ play a negative role in maximizing $Q(\beta,\beta')$.
This is done in the following lemma, stated in a more general setting and shown in Subsection \ref{subsection:random_perturbation}.
\begin{lem} \label{lem:random_perturbation}
Let $(p_n)_{n\geq 1}$ and $(q_n)_{n\geq 1}$ be nonincreasing deterministic sequences of nonnegative real numbers such that $\sum_{n\geq 1} p_n = 1$.
Let $(A_n)_{n\geq 1}$ be a sequence of i.i.d.\@ positive random variables.
We set 
\[
\tilde{p}_n \coloneqq \frac{A_n p_n}{\sum_{k\geq 1} A_k p_k}, 
\quad \Forall n\geq 1.
\]
Then, we have 
\begin{equation}
\Ec{\sum_{n\geq 1} \tilde{p}_n q_n} 
\leq \sum_{n\geq 1} p_n q_n. \label{eq:random_perturbation}
\end{equation}
Moreover, if $A_1$ is not almost surely constant, $(q_n)_{n\geq 1}$ is not constant and, for any $n \geq 1$, $p_n > 0$, then the inequality in \eqref{eq:random_perturbation} is strict.
\end{lem}
\begin{rem} In order to maximize the inner product $\sum_{n\geq 1} p_n q_n$, the sequences $(p_n)_{n\geq 1}$ and $(q_n)_{n\geq 1}$ have to be ordered in the same way, as in particular in the assumption of the lemma. 
The random perturbation on $(\tilde{p}_n)_{n\geq 1}$ may disturb this shared ordering and, hence, tends to reduce the inner product.
However, note that, without the expectation, the inequality in \eqref{eq:random_perturbation} can be wrong with positive probability. 
Therefore, apart from the case where $\beta$ or $\beta'$ is infinite (see Remark \ref{rem:beta_infinite}), it is not clear whether 
$Q(\beta,\beta')$ is stochastically dominated by $Q^\mathrm{REM}(\beta,\beta')$ or not.
\end{rem}
\begin{rem} One could also assume that $\sum_{n\geq 1} q_n = 1$ and perturb the sequence $(q_n)_{n\geq 1}$, by defining $\tilde{q}_n \coloneqq B_n q_n / \sum_{k\geq 1} B_k q_k$ with $(A_n,B_n)_{n\geq 1}$ a sequence of i.i.d.\@ random variables with values in $(0,\infty)^2$.
But then, we \textit{do not} necessarily have $\E[\sum_{n\geq 1} \tilde{p}_n \tilde{q}_n] \leq \sum_{n \geq 1} p_n q_n$.
For this reason, in order to prove Theorem \ref{theorem:comparaison_overlap_GFF}, we first have to use Lemma \ref{lem:changement_de_processus_ponctuel}, so that only one of the two sequences is perturbed (in comparison with the REM).
\end{rem}

Looking at Lemma \ref{lem:random_perturbation}, one can observe that, in order to get a strict inequality in Theorem~\ref{theorem:comparaison_overlap_GFF}, we still need the following lemma, which is proved in Subsection \ref{subsection:decorations_are_nontrivial}. 
\begin{lem} \label{lem:decoration_non_triviales}
The random variable $Y$ defined by \eqref{eq:def_Y} is not almost surely constant.
\end{lem}
This result is the only step in the proof of Theorem \ref{theorem:comparaison_overlap_GFF} which is specific to the DGFF, because it depends on the law $\nu$ of the decorations in the limit of the full extremal process. 
Now, with these three lemmas, we can proceed to the proof of Theorem \ref{theorem:comparaison_overlap_GFF}.
\begin{proof}[Proof of Theorem \ref{theorem:comparaison_overlap_GFF}]
Applying Lemma \ref{lem:changement_de_processus_ponctuel}, we get
\[
\Ec{Q(\beta,\beta')}
= \Ec{ \frac{\sum_{n\geq 1} \e^{\beta \xi_n} 
\e^{\beta' (\xi_n + Y_n)}}
{\left( \sum_{k\geq 1} \e^{\beta \xi_k} \right)
\left( \sum_{k\geq 1} \e^{\beta' (\xi_k + Y_k)} \right)}
}.
\]
We can assume that the atoms $(\xi_k)_{k\geq 1}$ are ranked in decreasing order.
Then, Lemma \ref{lem:random_perturbation} with
\[
p_n \coloneqq 
\frac{\e^{\beta' \xi_n}}{\sum_{k\geq 1} \e^{\beta' \xi_k}},
\quad
q_n \coloneqq 
\frac{\e^{\beta \xi_n}}{\sum_{k\geq 1} \e^{\beta \xi_k}},
\quad \text{and} \quad
A_n \coloneqq
\e^{\beta' Y_n},
\]
where $A_n$ is not almost surely constant by Lemma \ref{lem:decoration_non_triviales}, implies that, almost surely,
\[
\Ecsq{\frac{\sum_{n\geq 1} \e^{\beta \xi_n} 
\e^{\beta' (\xi_n + Y_n)}}
{\left( \sum_{k\geq 1} \e^{\beta \xi_k} \right)
\left( \sum_{k\geq 1} \e^{\beta' (\xi_k + Y_k)} \right)}
}{(\xi_k)_{k\geq 1}} 
< \frac{\sum_{n\geq 1} \e^{\beta \xi_n} 
\e^{\beta' \xi_n}}
{\left( \sum_{k\geq 1} \e^{\beta \xi_k} \right)
\left( \sum_{k\geq 1} \e^{\beta' \xi_k} \right)} = Q^\mathrm{REM}(\beta,\beta'),
\]
recalling the definition of $Q^\mathrm{REM}(\beta,\beta')$ from \eqref{eq:REM}.
Taking the expectation proves the result.
\end{proof}

\subsection{Change of point process}
\label{subsection:change_of_point_process}

\begin{proof}[Proof of Lemma \ref{lem:changement_de_processus_ponctuel}]
By \cite[Corollary 2.3]{kallenberg2017}, it is sufficient to prove that, for any $f \colon \R^2 \to \R_+$ continuous with compact support, 
\begin{align*}
\Ec{\exp \left( 
- \sum_{k\geq 1} f\left(\xi_k + X_{\beta,k},\xi_k + X_{\beta',k}\right) 
\right) }
& = \Ec{\exp \left( 
- \sum_{k\geq 1} f(\xi_k + c_\beta,\xi_k + c_\beta + Y_k)
\right) }.
\end{align*}
Starting from the left-hand side and integrating with respect to $(X_{\beta,k},X_{\beta',k})_{k\geq 1}$, we get
\begin{align} \label{ae}
\Ec{\exp \left( 
- \sum_{k\geq 1} f\left(\xi_k + X_{\beta,k},\xi_k + X_{\beta',k}\right) 
\right) }
& = \Ec{\exp \left( - \sum_{k\geq 1} \varphi(\xi_k) \right) },
\end{align}
where we set $\varphi(x) \coloneqq -\log \E[\e^{-f(x+X_{\beta},x+ X_{\beta'}) }]$ for any $x \in \R$.
Then, using the exponential formula for Poisson point processes, \eqref{ae} is equal to
\begin{align*}
\exp \left(
- \int_\R \diff x \, \e^{-\beta_c x} 
\left( 1 - \e^{-\varphi(x)} \right) \right) 
& = \exp \left(
- \Ec{ \int_\R \diff y \, \e^{-\beta_c (y-X_{\beta})} 
\left( 1 - \e^{-f(y,y+ X_{\beta'}-X_{\beta}) } \right)}\right),
\end{align*}
using Fubini's theorem and changing $x$ into $y-X_{\beta}$.
Then, applying again Fubini's theorem and recalling the definition of the random variable $Y$ in \eqref{eq:def_Y}, we get that \eqref{ae} is equal to
\begin{align} \label{ag}
\exp \left( - \int_\R \diff y \, \Ec{\e^{\beta_c X_{\beta}}} \e^{-\beta_c y} 
\left( 1 - \Ec{\e^{-f(y,y+ Y)}} \right) \right).
\end{align}
We proceed then in the same way, but backwards: setting $\psi(x) \coloneqq -\log \E[\e^{-f(x,x+Y) }]$ for any $x \in \R$ and changing $y$ into $z + c_\beta$, \eqref{ag} is equal to
\begin{align*}
\exp \left(
- \int_\R \diff z \, \e^{-\beta_c z} 
\left( 1 - \e^{-\psi(z + c_\beta)} \right) \right)
& = \Ec{\exp \left( - \sum_{k\geq 1} \psi(\xi_k + c_\beta) \right) } \\
& = \Ec{\exp \left( 
- \sum_{k\geq 1} f(\xi_k + c_\beta,\xi_k + c_\beta + Y_k)
\right) },
\end{align*}
which concludes the proof.
\end{proof}

\subsection{Random perturbation of a nonincreasing sequence}
\label{subsection:random_perturbation}

\begin{proof}[Proof of Lemma \ref{lem:random_perturbation}]
If $\sum_{n\geq1} A_n p_n = \infty$ a.s.\@, then $\tilde{p}_n = 0$ a.s.\@ for any $n \geq 1$ and therefore \eqref{eq:random_perturbation} holds clearly and the inequality is strict as soon as $(q_n)_{n\geq 1}$ is not constant equal to zero.
Otherwise, we have $\sum_{n\geq1} A_n p_n < \infty$ a.s.\@ by Kolmogorov's zero-one law and we work under this assumption until the end of the proof. Using also that $\sum_{n\geq1} p_n = 1$ and $(q_n)_{n\geq 1}$ is bounded, it follows that all forthcoming sums converge absolutely almost surely.
Note also that if $(q_n)_{n\geq 1}$ is constant, then \eqref{eq:random_perturbation} is clear, hence we assume that it is not.
We introduce 
\[
S \coloneqq \sum_{n\geq 1} (p_n-\tilde{p}_n) q_n
\]
and our aim is to prove that $\E[S] \geq 0$.
Using that $\sum_{k\geq 1} p_k =1$, we have 
\[
S = \sum_{n\geq 1} 
\frac{\left(\sum_{k\geq 1} A_k p_k \right) p_n
- \left(\sum_{k\geq 1} p_k \right) A_n p_n}
{\sum_{j\geq 1} A_j p_j} q_n 
= \frac{1}{\sum_{j\geq 1} A_j p_j}
\sum_{n,k \geq 1} (A_k - A_n) p_k p_n q_n.
\]
Then, switching the role of $k$ and $n$ in one of the sums, we get 
\[
\sum_{n,k \geq 1} (A_k - A_n) p_k p_n q_n
= \sum_{n,k \geq 1} A_k p_k p_n q_n
- \sum_{n,k \geq 1} A_n p_k p_n q_n
= \sum_{n,k \geq 1} A_n p_n p_k (q_k - q_n).
\]
Therefore, applying dominated convergence theorem, we have 
\[
\E[S] 
= \sum_{n,k \geq 1} \Ec{\frac{A_n}{\sum_{j\geq 1} A_j p_j} } 
	p_n p_k (q_k - q_n)
= \sum_{n \geq 1} y_n p_n x_n,
\]
where we set, for any $n \geq 1$,
\[
x_n \coloneqq \sum_{k \geq 1} p_k (q_k - q_n)
\quad \text{and} \quad 
y_n \coloneqq  \Ec{\frac{A_n}{\sum_{j\geq 1} A_j p_j}}.
\]
Note that $x_n = (\sum_{k\geq 1} q_k p_k) - q_n$ is nondecreasing in $n$, with $x_1 \leq 0$ and $\lim_{n\to\infty} x_n >0$ (because $\sum_{k\geq 1} p_k = 1$ and $(q_n)_{n\geq 1}$ is nonincreasing and not constant), and therefore 
there exists $n_0 \geq 1$ such that $x_n \leq 0$ for $n \leq n_0$ and $x_n > 0$ for $n > n_0$.
Moreover, we have $\sum_{n\geq 1} p_n x_n = 0$, where the $n_0$ first terms are nonpositive and the other ones are nonnegative.
Therefore, in order to prove that $\E[S] = \sum_{n \geq 1} y_n p_n x_n$ is nonnegative, it is sufficient to prove that $(y_n)_{n\geq 1}$ is a nondecreasing sequence.
Indeed, it follows then that
\[
\sum_{n\geq 1} y_n p_n x_n 
\geq \sum_{n\geq 1} y_{n_0} p_n x_n 
= 0,
\]
by distinguishing the case $n \leq n_0$ and the case $n > n_0$.
For any $n>m \geq 1$, our aim is now to prove that $y_n \geq y_m$. 
Setting $B \coloneqq \sum_{j\notin \{m,n\}} A_j p_j$, we have
\[
y_n - y_m 
= \Ec{\frac{A_n-A_m}{A_n p_n + A_m p_m + B}}.
\]
We distinguish two cases.
On the event $\{ A_n \geq A_m \}$, since $p_n \leq p_m$, we have $A_n p_n + A_m p_m \leq (A_n+A_m)(p_n+p_m)/2$ and so
\begin{equation}
\frac{A_n-A_m}{A_n p_n + A_m p_m + B}
\geq \frac{A_n-A_m}{\frac{(A_n + A_m)(p_n + p_m)}{2} + B}. \label{ab}
\end{equation}
On the other hand, on the event $\{ A_n < A_m \}$, we have $A_n p_n + A_m p_m \geq (A_n+A_m)(p_n+p_m)/2$ and \eqref{ab} holds too.
Therefore, we get
\begin{equation}
y_n - y_m 
\geq \Ec{\frac{2(A_n-A_m)}{(A_n + A_m)(p_n + p_m) + 2B}} = 0,
\label{ac}
\end{equation}
because $A_n$ and $A_m$ have the same distribution and play a symmetric role in the last expectation (and $B$ is independent of $(A_n,A_m)$). 
It proves that $(y_n)_{n\geq 1}$ is a nondecreasing sequence and therefore that $\Ec{S} \geq 0$.
Now, we assume in addition that $A_1$ is not almost surely constant and, for each $n \geq 1$, $p_n > 0$. Our aim is to prove that
\begin{equation} \label{ad}
\E[S] = \sum_{n \geq 1} y_n p_n x_n > 0.
\end{equation}
First note that, for any $n>m \geq 1$, if $p_n < p_m$ then $y_n > y_m$. 
Indeed, on the event $\{A_n < A_m \}$, we have in this case $A_n p_n + A_m p_m > (A_n+A_m)(p_n+p_m)/2$ and, therefore,
\[
\frac{A_n-A_m}{A_n p_n + A_m p_m + B}
> \frac{A_n-A_m}{\frac{(A_n + A_m)(p_n + p_m)}{2} + B}. 
\]
But the event $\{A_n < A_m \}$ has positive probability (because $A_n$ and $A_m$ are independent with the same non-constant distribution) and, thus, we get a strict inequality in \eqref{ac}.
We can now prove \eqref{ad}.
Recall that $\sum_{n\geq 1} p_n x_n = 0$, where the $n_0$ first terms are nonpositive and the other ones are positive (because now $p_n >0$).
Since $p_n \downarrow 0$ and $p_{n_0} > 0$, there exists $n_1 > n_0$ such that $p_{n_1} < p_{n_0}$ and therefore $y_{n_1} > y_{n_0}$.
Then, using that $p_{n_1} x_{n_1} > 0$, we get the following \textit{strict} inequality
\[
\sum_{n\geq 1} y_n p_n x_n 
> y_{n_0} \sum_{n\geq 1} p_n x_n 
= 0,
\]
which concludes the proof.
\end{proof}

\subsection{Decorations play a nontrivial role}
\label{subsection:decorations_are_nontrivial}

In this subsection we prove Lemma \ref{lem:decoration_non_triviales}, which states that the random variable $Y$, whose distribution depends on the decoration law $\nu$, is not almost surely constant. The contrary would have been surprising, but, in order to prove it properly, we first need to show a basic property of the law $\nu$ in Lemma \ref{lem:decoration_with_density} below.

Let us start by recalling some results from \cite{biskuplouidor2018} concerning the decoration law $\nu$ for the DGFF.
Let $\nu^0$ denote the law of the DGFF in $\Z^2 \setminus \{0\}$.
Its covariance is given by $\Cov_{\nu^0} (\phi_x,\phi_y) = \mathfrak{a}(x) + \mathfrak{a}(y) - \mathfrak{a}(x-y)$, where $\mathfrak{a} \colon \Z^2 \to \R$ is the potential kernel of the simple symmetric random walk started from zero (see Appendix \ref{appendix} for more details).
In this subsection, we denote by $\E_\nu$ and $\E_{\nu^0}$ the expectations under which $\phi = (\phi_x)_{x \in \Z^2}$ has respectively law $\nu$ and $\nu^0$.
Let $\cC_b^\mathrm{loc}(\R^{\Z^2})$ denote the set of continuous bounded functions on $\R^{\Z^2}$ that depend only on a finite number of coordinates in $\Z^2$.
Then, Biskup and Louidor \cite[Theorem 2.3 and Proposition 5.8]{biskuplouidor2018} proved the following description for the decoration law: for any $f \in \cC_b^\mathrm{loc}(\R^{\Z^2})$,
\begin{align} \label{eq:description_decoration_law}
\E_{\nu} \left[ f(\phi) \right]
= \lim_{r\to\infty} \E_{\nu^0} \left[
f (\psi)
\mathrel{}\middle|\mathrel{}
\Forall x \in \Lambda_r(0), \psi_x \geq 0
\right],
\qquad \text{where }
\psi \coloneqq \phi + \frac{2}{\sqrt{g}} \mathfrak{a},
\end{align}
recalling that $\Lambda_r(0) = \{x \in \Z^2 : \norme{x} \leq r \}$.
The main tool for the proof of Lemma \ref{lem:decoration_non_triviales} is the following lemma, that describes the conditional law of one coordinate of the decoration given the other coordinates. 
\begin{lem} \label{lem:decoration_with_density}
For any $y \in \Z^2\setminus\{0\}$, under $\nu$, the conditional law of $\phi_y$ given $(\phi_x)_{x \neq y}$ is a normal law with mean $\frac{1}{4} \sum_{x \sim y} \phi_x$ and variance $1$ conditioned to be nonnegative.
\end{lem}
\begin{proof}
Let $f \in \cC_b^\mathrm{loc}(\R^{\Z^2})$ that does not depend on the $y$-coordinate. 
Let $h \in \cC_b(\R)$.
Applying the domain Markov property to the DGFF $\phi$ with law $\nu^0$, 
conditionnally on the $\phi_x$ for $x \neq y$, 
$\phi_y$ has the same law as $Z + \frac{1}{4} \sum_{x \sim y} \phi_x$, 
where $Z$ is independent of $(\phi_x)_{x \neq y}$ with law $\cN(0,1)$.
On the other hand, by \cite[Proposition 4.4.2]{lawlerlimic2010}, the potential kernel $\mathfrak{a}$ is discrete harmonic on $\Z^2\setminus\{0\}$, so we have $\mathfrak{a}(y) = \frac{1}{4} \sum_{x \sim y} \mathfrak{a}(x)$.
Therefore, we get
\begin{align}
& \E_{\nu^0} \Biggl[
f(\psi) h(\psi_y)
\prod_{x \in \Lambda_r(0)} \1_{\{\psi_x \geq 0\}}
\Biggr] \nonumber \\
& = \E_{\nu^0} \Biggl[
f(\psi) 
\Biggl(
	\prod_{x \in \Lambda_r(0) \setminus \{y\}} \1_{\{\psi_x \geq 0\}}
\Biggr)
\int_\R
h \left( z + \frac{1}{4} \sum_{x \sim y} \psi_x \right)
\1_{\{z + \frac{1}{4} \sum_{x \sim y} \psi_x \geq 0\}}
\frac{\e^{-z^2/2}}{\sqrt{2\pi}} \diff z
\Biggr]. \label{ca}
\end{align}
Now, we define a function $H \colon \R \to \R$ by 
\begin{align*}
H(t)
& \coloneqq 
\frac{\int_\R h(z+t) \1_{\{z+t \geq 0\}}
\frac{\e^{-z^2/2}}{\sqrt{2\pi}} \diff z }
{\int_\R \1_{\{z+t \geq 0\}}
\frac{\e^{-z^2/2}}{\sqrt{2\pi}} \diff z},
\end{align*}
which is a continuous bounded function.
Then, \eqref{ca} is equal to
\begin{align}
& \E_{\nu^0} \Biggl[
f(\psi)
\Biggl(
	\prod_{x \in \Lambda_r(0) \setminus \{y\}} \1_{\{\psi_x \geq 0\}}
\Biggr)
H \left( \frac{1}{4} \sum_{x \sim y} \psi_x \right)
\int_\R
\1_{\{z + \frac{1}{4} \sum_{x \sim y} \psi_x \geq 0\}}
\frac{\e^{-z^2/2}}{\sqrt{2\pi}} \diff z
\Biggr] \nonumber \\
& = \E_{\nu^0} \Biggl[
f(\psi)
H \left( \frac{1}{4} \sum_{x \sim y} \psi_x \right)
\prod_{x \in \Lambda_r(0)}
\1_{\{\psi_x \geq 0\}}
\Biggr], \label{cb}
\end{align}
using the Markov domain property as before.
Applying \eqref{eq:description_decoration_law} to the left-hand side of \eqref{ca} and to the right-hand side of \eqref{cb}, we get
$\E_{\nu}[ f(\phi) h \left( \phi_y \right)] 
= \E_{\nu}[ f(\phi) H(\frac{1}{4} \sum_{x \sim y} \phi_x)]$.
Recalling the definition of function $H$ and since this equality holds for any $f \in \cC_b^\mathrm{loc}(\R^{\Z^2})$ not depending on the $y$-coordinate and any $h \in \cC_b(\R)$, it proves the claimed result.
\end{proof}
\begin{proof}[Proof of Lemma \ref{lem:decoration_non_triviales}] 
By contradiction, assume that $Y$ is almost surely constant.
Then, $X_{\beta'} - X_\beta$ is a.s.\@ constant and, taking the exponential, there exists some constant $c > 0$ such that
\begin{align*}
\sum_{x \in \Z^2} \e^{-\beta \phi_x} 
= c \left( \sum_{x \in \Z^2} \e^{-\beta' \phi_x} \right)^{\beta/\beta'},
\quad \nu\text{-a.s.}
\end{align*}
Fix some $y \in \Z^2 \setminus \{0\}$.
We have, $\nu$-a.s.,
\begin{align*}
\E_{\nu} \left[ 
\abs{ \e^{-\beta \phi_y} + \sum_{x \neq y} \e^{-\beta \phi_x} 
- c \Biggl( \e^{-\beta \phi_y} + \sum_{x \neq y}  \e^{-\beta' \phi_x} \Biggr)^{\beta/\beta'} }
\mathrel{}\middle|\mathrel{}
(\phi_x)_{x \neq y} 
\right]
= 0,
\quad \nu\text{-a.s.}
\end{align*}
But, given $(\phi_x)_{x \neq y}$, $\phi_y$ has a distribution with positive density on $(0,\infty)$ by Lemma \ref{lem:decoration_with_density}. Recalling that $\sum_{x \neq y} \e^{-\beta \phi_x} \geq 1$, it shows that there exist some $s,s' \geq 1$ such that, 
\begin{align*}
\e^{-\beta z} + s - c \left( \e^{-\beta z} + s' \right)^{\beta/\beta'} = 0,
\quad \text{for almost every } z \in (0,\infty).
\end{align*}
This cannot be true for $\beta \neq \beta'$ and, therefore, one gets a contradiction.
\end{proof}

\appendix

\section{Some useful results}
\label{appendix}

In this appendix, we state some known results used throughout the paper.
The first result is the well-known {\it Gaussian integration by part}, whose proof can be found in \cite[Equation (A.17)]{talagrand2011}.
\begin{lem}
\label{lem:Gaussian-IbP} 
Let $(X,Z_1,\dots , Z_d)$ be a centered Gaussian random vector. Then, for any $\cC^1$ function $F \colon \R^d \rightarrow \R$ of moderate growth at infinity, we have 
\begin{equation*}
\E \left[ X F(Z_1,\dots,Z_d) \right] 
= \sum_{i=1}^d  \E \left[ X Z_i \right] \cdot
\E \left[ \frac{\partial F}{\partial z_i}(Z_1,\dots,Z_d) \right].
\end{equation*}
\end{lem}
Now we define formally the Green function $G_N$ and state some of its properties.
Let $(S_n)_{n\geq 0}$ be a simple symmetric random walk on $\Z^2$ starting from $x$ under $\P_x$ and $\tau_N \coloneqq \inf \{ n \geq 0 : S_n \notin D_N \}$ be the first exit time of the walk from the set $D_N$.
Then, the Green function is defined, for $x,y \in \Z^2$, by
\begin{equation}
\label{eq:Greenfunction}
G_N(x,y) = \E_x \left[ \sum_{k=0}^{\tau_N-1} \1_{\{S_k = y\}} \right].
\end{equation}
The Green function is related to the potential kernel of the simple symmetric random walk started from zero, which is a function  $\mathfrak{a} \colon \Z^2 \to \R$ defined by
\[
\mathfrak{a}(x) \coloneqq \sum_{n\geq 0} (\P(S_n = 0) - \P(S_n = x))
\]
and can also been written explicitely as a double integral, see \cite[Proposition 4.4.3]{lawlerlimic2010}.
Indeed, by \cite[Proposition 4.6.2]{lawlerlimic2010} and since $D_N$ is finite, the following relation between these functions holds for any $x,y \in \Z^2$,
\begin{equation} \label{eq:relation_G_N_and_a}
G_N(x,y) = \E_x \left[ \mathfrak{a}(S_{\tau_N} - y) \right] - \mathfrak{a}(x-y).
\end{equation}
Furthermore, Theorem 4.4.4 of \cite{lawlerlimic2010} gives the asymptotic behavior of function $\mathfrak{a}$: as $\norme{x} \to \infty$, we have
\begin{equation} \label{eq:asymptotics_a}
\mathfrak{a} (x) = \frac{2}{\pi} \log \norme{x} + \frac{2 \gamma + \log 8}{\pi} + \grandO{\norme{x}^{-2}},
\end{equation}
where $\gamma$ is Euler's constant.
From these tools, we can easily deduce the following result, relating the overlap between two points to their euclidean distance.
\begin{lem} \label{lemma:overlap-distance}  
\begin{enumerate}
\item There exists a constant $c_1 > 0$ such that, for any $N \geq 1$, $0<\alpha<1$ and any $x,y \in D_N$,
\[
q_N (x,y) \geq \alpha \quad \Rightarrow \quad 
\norme{x-y} \le c_1 N^{1-\alpha}.
\]
\item There exists a constant $c_2 > 0$ such that, for any $N \geq 1$, $0<\delta<\alpha<1$ and any $x,y \in D_{N,\delta} \coloneqq \{ x \in D_N : d(x,D_N^c) > N^{1-\delta} \}$,
\[
\norme{x-y} \le c_2 N^{1-\alpha} \quad \Rightarrow \quad 
q_N (x,y) \geq \alpha-\delta.
\]
\end{enumerate}
\end{lem}

\begin{proof}
Recall that $q_N(x,y) = G_N(x,y) / \max_{z \in D_N} G_N(z,z)$.
Combining \eqref{eq:relation_G_N_and_a} and \eqref{eq:asymptotics_a}, first note that $\max_{z \in D_N} G_N(z,z) = \frac{2}{\pi} \log N + O_N(1)$ as $N \to \infty$.
Using again \eqref{eq:relation_G_N_and_a} and \eqref{eq:asymptotics_a}, we get the following upper bound, uniformly in $x,y \in D_N$,
\begin{equation} \label{eq:upper_bound_overlap}
q_N(x,y) \leq 1-\frac{\log \Vert x-y \Vert}{\log N} +O_N((\log N)^{-1}),
\end{equation}
and the following lower bound, uniformly in $x,y \in D_{N,\delta}$,
\begin{equation} \label{eq:lower_bound_overlap}
q_N(x,y) \geq 1- \delta - \frac{\log \Vert x-y \Vert}{\log N} +O_N((\log N)^{-1}). 
\end{equation}
Part (i) of the lemma follows from \eqref{eq:upper_bound_overlap} and Part (ii) from \eqref{eq:lower_bound_overlap}.
\end{proof}
In the previous lemma, we saw that the overlap is properly related to the distance if the points of interest are not too close from the boundary of $D_N$.
The following result shows that, for any inverse temperature $\beta>0$, the Gibbs measure of $D_{N,\delta}^c \coloneqq D_N \setminus D_{N,\delta}$ is negligible.
\begin{lem} \label{lemma:complement-negligible}
Let $\delta > 0$ and recall that $D_{N,\delta} \coloneqq \{ x \in D_N : d(x,D_N^c) > N^{1-\delta} \}$.
For any $\beta >0$,
 \[
\cG_{\beta,N} (D_{N,\delta}^c)
\xrightarrow[N \to\infty]{} 0, \quad \textrm{in } L^1.
\]
\end{lem}

\begin{proof}
This result is shown in \cite[Lemma 3.1]{arguinzindy2015} in the case of the square lattice $V_N = (0,N)^2 \cap \Z^2$, but the proof works also in our more general framework because we also have $\# D_{N,\delta}^c = O(N^{2-\delta})$ under our assumptions.
Indeed, we assumed that $\partial D$ has a finite number of connected components, each of which is a $\cC^1$ path: therefore, there exists $C > 0$ such that, for any $\varepsilon > 0$, $\partial D$ can be covered by at most $C \varepsilon^{-1}$ balls of radius $\varepsilon$.
Fix some $N \geq 1$, there exist $M \leq CN^\delta /2$ and $z_i,\dots,z_M$ such that $\partial D \subset \bigcup_{i=1}^M B(z_i,2N^{-\delta})$.
Then, recalling the definition of $D_N$ and $D_{N,\delta}$, note that
\[
\frac{D_{N,\delta}^c}{N} 
\subset \left\{ z \in D : d(z,\partial D) \leq 2N^{-\delta} \right\} \cap \frac{\Z^2}{N}
\subset \bigcup_{i=1}^M \left( B(z_i,4N^{-\delta}) \cap \frac{\Z^2}{N} \right).
\]
It follows that $\# D_{N,\delta}^c \leq M 4 (4 N^{1-\delta})^2 \leq C' N^{2-\delta}$, which is the announced result.
\end{proof}

\bigskip

\section*{Acknowledgements}

We thank Antonio Auffinger for suggesting the question of chaos in temperature for the two-dimensional discrete Gaussian free field and Bernard Derrida for very stimulating discussions. The authors also wish to thank the referees for useful comments that improved the presentation of the paper.

\addcontentsline{toc}{section}{References}

\bibliographystyle{abbrv}


\end{document}